\documentclass[12pt,a4paper]{article}
\usepackage[catalan,spanish,english]{babel}
\usepackage[psamsfonts]{amssymb}
\usepackage{cmmib57}
\usepackage{exscale}
\usepackage{amsmath}
\usepackage{amscd}
\usepackage{latexsym}
\usepackage{graphicx}
\newtheorem{thm}{Theorem}[section]
\newtheorem{cor}[thm]{Corollary}

\newtheorem{prop}[thm]{Proposition}

\newtheorem{rem}{Remark}[section]
\newcommand{\qed}{\quad{$\square$}}

\newcommand{\Ind}{{\rm l}\hskip -0.23truecm 1}
\newcommand{\MultipleIntegral}{\mathbf{I}}
\newcommand{\fnint}{\MultipleIntegral^{(n),H}}

\newcommand{\bh}{B^H}
\newcommand{\kh}{K_H}
\newcommand{\kht}{K_{H,t}}

\newcommand{\khst}{K_{H,s,t}}
\newcommand{\rkhsh}{\mathcal{H}^H}

\newcommand{\IntervalT}{\left[0,T\right]}
\newcommand{\Intervalt}{\left[0,t\right]}
\newcommand{\Partialr}{\frac{\partial \kh}{\partial r}}
\newcommand{\Partialt}{\frac{\partial \kh}{\partial t}}
\newcommand{\RegularSpace}{\left|\rkhsh\right|}
\newcommand{\SingularSpace}{ \mathcal{H}_{\kh} }
\newcommand{\StepFunctions}{\mathcal{E}}
\newcommand{\SobolevSpace}{\mathbb{D}}
\newcommand{\nLambda}{\Lambda^{(n)}}

\newcommand{\nRegularSpace}{L^{\frac{1}{H}}(\nLambda_T)}
\newcommand{\nRSHolder}{L^q(\nLambda_T)}
\newcommand{\nSingularSpace}{\mathbb{H}^{\lambda}(\nLambda_T )}

\newcommand{\RegularSpaceIntegrand}{L^{\frac{1}{H}}\left(\IntervalT;\SobolevSpace^{1,p}\right)}
\newcommand{\SingularSpaceIntegrand}{\mathcal{C}^{\lambda}\left(\IntervalT;\SobolevSpace^{1,p}\right)}

\newcommand{\ep}{\varepsilon}
\newcommand{\sep}{\sqrt{\varepsilon}}

\begin{document}
\begin{titlepage}
\null
\vspace{2cm}
\begin{center}
{\Large\bf A Large Deviation Principle in H\"older Norm\\
for Multiple Fractional Integrals \\[2mm]}

\vspace{7mm}

\begin{tabular}{l@{\hspace{10mm}}l@{\hspace{10mm}}l}
{\sc Marta Sanz-Sol\'e}$\,^{(\ast)}$ &and &{\sc Iv\'an Torrecilla-Tarantino}$\,^{(\ast)}$\\
{\small marta.sanz@ub.edu }         &&{\small itorrecilla@ub.edu}\\
\end{tabular}
\begin{center}
{\small Facultat de Matem\`atiques}\\
{\small Universitat de Barcelona } \\
{\small Gran Via 585} \\
{\small 08007 Barcelona, Spain} \\
\end{center}
\end{center}

\vspace{1.5cm}

\noindent{\bf Abstract:} For a fractional Brownian motion $B^H$
with Hurst parameter
$H\in]\frac{1}{4},\frac{1}{2}[\cup]\frac{1}{2},1[$, multiple
indefinite integrals on a simplex are constructed and the
regularity of their sample paths are studied. Then, it is proved
that the family of probability laws of the  processes obtained by
replacing $B^H$ by $\varepsilon^{\frac{1}{2}} B^H$ satisfies a
large deviation principle  in H\"older norm. The definition of the
multiple integrals relies upon a representation of the fractional
Brownian motion in terms of a stochastic integral with respect to
a standard Brownian motion. For the large deviation principle, the
abstract general setting in \cite{Le90} is used.
\medskip



\vspace{2 cm}

\noindent

\footnotesize

{\begin{itemize} \item[$^{(\ast)}$] Supported by the grant MTM 2006-01351 from the \textit{Direcci\'on General de
Investigaci\'on, Ministerio de Educaci\'on y Ciencia.}
\end{itemize}}
\end{titlepage}

\newpage

\section{Introduction}
\label{s0}
In this paper, we consider stochastic processes $X = \{X_t, t\in[0,T]\}$ given by indefinite multiple integrals on the $n$--dimensional simplex
$\{(\theta_1,\dots,\theta_n)\in \mathbb{R}_+^n: 0\le\theta_1\le\dots\le \theta_n\le t\}$ with respect to a fractional Brownian motion (fBm).
The integrands are allowed to depend also on the parameter $t$. Under suitable assumptions on the integrands, depending whether the Hurst parameter $H$ belongs to
$]\frac{1}{2},1[$ or $]\frac{1}{4}, \frac{1}{2}[$, we prove H\"older continuity of the sample paths, a.s..
Then we establish  a large deviation principle (LDP) in H\"older norm for the family of laws of $\varepsilon^{\frac{n}{2}} X$.

For the standard Brownian motion (sBm), a similar question has been addressed in \cite{MWNPA92}. The authors consider different assumptions on the  integrands
ensuring a.s. continuity of the  sample paths of the integrals;  then they prove large deviations principles in the space of continuous functions
endowed with the supremum norm.

Geometric rough paths based on processes with $\gamma$--H\"older
continuous sample paths give rise to random vectors whose
components are multiple Stratonovich integrals up to  order
$[\frac{1}{\gamma}]$, where $[\cdot]$ denotes the integer value. A
large deviation principle for the rough path lying above the
standard Brownian motion has been proved in \cite{LQZ}. The norm
under consideration is the $p$--variation norm used in the rough
path analysis (see for instance \cite{L-Q}). The higher order of
the multiple stochastic integrals involved is in this example
$n=2$. For the fractional Brownian motion with Hurst parameter
$H\in]\frac{1}{4}, \frac{1}{2}[\cup]\frac{1}{2},1[$, a similar
result has been proved in \cite{millet-ss}. We notice that the non
trivial part of it corresponds to the values $H\in]\frac{1}{4},
\frac{1}{2}[$ and that one needs to deal with multiple stochastic
integrals up to order $n=3$.

This paper is motivated mainly by  \cite{MWNPA92} and \cite{millet-ss} in the following sense: As in \cite{MWNPA92}, we want to consider multiple indefinite integrals of any order
and on the other hand, we wish to deal with sharper norms, like H\"older norm, and with the fBm.

The main body of the paper is devoted to the construction of the indefinite multiple integral on a simplex with respect to the fBm, and the study of
its sample paths. The corresponding results are gathered in Section \ref{s2}.
Starting from the results and ideas in \cite{AMN01} for $n=1$, by means of a recursive argument, we are able to give a meaning to the
multiple integral with respect to the fBm as a multiple integral with respect to the sBm. For this, we identify the kernels corresponding to increments
in time of such integrals. With suitable assumptions, we prove that these kernels define continuous operators on the space of the integrands
taking values on spaces of H\"older-continuous functionals (see (\ref{2.1.0.1}) and (\ref{2.2.4})). By means of the hypercontractivity property
of Gaussian chaos, the H\"older continuity is transferred to the sample paths of the integrals.

We should mention the fractional calculus approach (see for
instance \cite{SKM93}) to multiple definite integrals with respect
to the fBm given in  \cite{PAT02}, and  to  indefinite integrals
of progressively measurable processes with respect to the fBm
-including sample path properties- in \cite{D1}. In contrast, as
we have mentioned before, our approach follows  \cite{AMN01} (see
also \cite{DU198} and \cite{D2}). It is based on anticipating
integrals of Skorohod type; thus, on techniques from Malliavin
calculus.

Once we have identified the functional spaces where our fBm
functionals live, we can study what LDP they do satisfy. In
\cite{Le90},  a LDP for random vectors in a Banach-valued
homogeneous Wiener chaos of any order $n$ is established. The
elegant proof relies upon isoperimetric methods. This provides the
suitable framework for our study. In fact, in Section \ref{s3} we
first identify the abstract Wiener space associated with the fBm
as a Gaussian process with H\"older continuous paths. Then we
notice that the space of $\gamma$--H\"older continuous functions
can be embedded in a separable Banach space (see the first part of
Section \ref{s3} for some details and references). With this and
the results of Section \ref{s2}, we see  that the indefinite
multiple integrals with respect to the fBm are Banach-valued
random vectors in a Wiener chaos. Therefore, the results of
\cite{Le90} can be applied. A similar approach could be used for
the sBm to obtain the LDP stated in \cite{MWNPA92} and very likely
with sharper norms.

\section{Preliminaries and notation}
\label{s1}

We start the article with this section devoted to fix the notation and recall some known facts that will be intensively used
throughout the paper, refereeing  to \cite{Nu06} and the references herein for additional details.

Let $\bh=\left\{\bh_t,\;t \in \IntervalT \right\}$ be a fBm with Hurst parameter $H\in ]0,1[$.  The process $\bh$ can be
represented in terms of a stochastic integral with respect to a sBm $W=\left\{W_t,\;t \in \IntervalT\right\}$ as follows:
\begin{equation}
\bh_t=\int_0^t\kh(t,\theta)dW_\theta,\label{1.1}
\end{equation}
where $d W_\theta$ denotes the It\^o differential and
\begin{align}
\kh(t,\theta):=c_H\left\{\left(t-\theta\right)^{H-\frac{1}{2}}+
\left(\frac{1}{2}-H\right)\int_\theta^t\left(u-\theta\right)^{H-\frac{3}{2}}
\left(1-\left(\frac{\theta}{u}\right)^{\frac{1}{2}-H}\right)du\right\},\label{1.2}
\end{align}
 and $c_H$ is some positive constant depending on $H$.
Then
\begin{align}
\Partialt(t,\theta)=c_H\left(H-\frac{1}{2}\right)\left(\frac{\theta}{t}\right)^{\frac{1}{2}-H}\left(t-\theta\right)^{H-\frac{3}{2}}.
\label{1.3}
\end{align}
Thus, for $H\in]0,\frac{1}{2}[$, the derivative $\Partialt(t,\theta)$ is negative and moreover,
\begin{equation}
\label{1.20}
\left\vert \Partialt(t,\theta)\right\vert\le C_H \vert t-\theta\vert^{H-\frac{3}{2}}.
\end{equation}

Notice that for  $H>\frac{1}{2}$  the kernel
$\kh\left(t,\theta\right)$ is regular and for $H<\frac{1}{2}$ it
is singular.

For any $h\in L^2([0,T])$ we define the operator $K_H$ by
\begin{equation*}
\left(K_H h\right)(t) = \int_0^t K_H(t,\theta) h(\theta) d \theta.
\end{equation*}

Let $\StepFunctions$ be the set of step functions on
$\IntervalT$. We define the $L^2\left(\IntervalT\right)$--valued linear operator $\kh^\ast$ on
$\StepFunctions$ by
\begin{equation}
\label{1.3.1}
\left(\kh^\ast\varphi\right)\left(\theta\right)=\varphi\left(\theta\right)\kh\left(T,\theta\right)
+\int_\theta^T\left[\varphi\left(r\right)-\varphi\left(\theta\right)\right] K_H(dr, \theta).
\end{equation}

The operator $\kh^\ast$ is the adjoint of $\kh$ in the following
sense:

\medskip

For any function $\varphi \in \StepFunctions$ and $h \in
L^2\left(\IntervalT\right)$, one has
\begin{equation}
\int_0^T\left(\kh^\ast\varphi\right)\left(r\right)h\left(r\right)dr=\int_0^T\varphi
\left(r\right)\left(\kh h\right)\left(dr\right). \label{1.4}
\end{equation}

Replacing $h(s)ds$ by $d W_s$, with the same proof of (\ref{1.4})
it can be checked that for any $\varphi \in \StepFunctions$ the
element $\bh\left(\varphi\right):= \int_0^T \varphi(\theta) d
B^H_\theta$ of the first Gaussian chaos associated with the fBm
can be written as
\begin{equation*}
\bh\left(\varphi\right)=\int_0^T\left(\kh^\ast\varphi\right)\left(\theta\right)d
W_\theta.
\end{equation*}

Let $t\in [0,T]$. Then
\begin{align}
\label{1.6}
\left(\kht^\ast\varphi\right)(\theta):&=\left(\kh^{\ast}(\varphi\;\Ind_{\left[0,t\right]})\right)(\theta)\nonumber\\
&=\varphi\left(\theta\right)\kh\left(t,\theta\right)
+\int_\theta^t\left[\varphi\left(r\right)-\varphi\left(\theta\right)\right]
K_H(dr,\theta).
\end{align}
Notice that $K^{\ast}_{H,T}= K^{\ast}_{H}$.

Thus,
\begin{equation*}
\bh\left(\varphi\;\Ind_{\left[0,t\right]}\right)=\int_0^t
\varphi(\theta) dB^H_\theta
=\int_0^t\left(\kht^\ast\varphi\right)\left(\theta\right)dW_\theta.
\end{equation*}
For $H\in]\frac{1}{2}, 1[$ and $\varphi \in \StepFunctions$, the kernel $\kht^\ast$ has the simple expression
\begin{equation}
\label{1.7.2} \left(\kht^\ast
\varphi\right)\left(\theta\right)=\int_\theta^t \varphi\left(r\right) K_H(dr,\theta).
\end{equation}

For this same range of $H$, denote by $\RegularSpace$ the linear space consisting of measurable functions
$\varphi$ defined on $\IntervalT$ such that
\begin{align*}
\left\|\varphi\right\|^2_{\RegularSpace}&:=
\int_0^T \left(\int_\theta^T\left|\varphi_r\right| K_H(dr,\theta)\right)^2d\theta\\
&=\alpha_H\int_0^T
\int_0^T\left|\varphi_r\right|\left|\varphi_\xi\right|\left|r-\xi\right|^{2H-2}drd\xi<\infty,
\end{align*}
where $\alpha_H=H\left(2H-1\right)$. The space $\RegularSpace$ endowed with the norm $\left\|\cdot\right\|_{\RegularSpace}$, is a Banach
space.
H\"{o}lder's inequality with exponent
$q=\frac{1}{H}$ and the Hardy-Littlewood-Sobolev inequality (see for instance, \cite{Stein}, page 354) yields
\begin{equation}
\left\|\varphi\right\|_{\RegularSpace}\leq b_H
\left\|\varphi\right\|_{L^{\frac{1}{H}}\left(\IntervalT\right)}.
\label{1.8}
\end{equation}
For $H\in]0,\frac{1}{2}[$, we introduce the seminorm on
$\StepFunctions$
\begin{equation*}
\left\|\varphi\right\|^2_{\kh}:=\int_0^T\varphi_\theta^2\kh\left(T,\theta\right)^2d\theta+
\int_0^T\left(\int_\theta^T\left|\varphi_r-\varphi_\theta\right|\left\vert K_H(dr,\theta\right\vert) \right)^2d\theta.
\end{equation*}
By $\SingularSpace$, we denote the completion of $\StepFunctions$
with respect to this seminorm. It consists of functions $\varphi$
defined on $\IntervalT$ such that $\left\|\varphi\right\|_{\kh}^2<
\infty$.

For any $t\in[0,T]$, set
$$\nLambda_t= \{(\theta_1,\dots,\theta_n)\in \mathbb{R}_+^n: 0\le\theta_1\le\dots\le \theta_n\le t\}.$$
Throughout the paper, we denote by $\nSingularSpace$,
$\lambda\in]0,1[$, the space of  $\lambda$--H\"older continuous
functions on the $k$--cubes contained on $\nLambda_t$, $1\leq k
\leq n$, endowed with the norm
\begin{align*}
&\left\|h\right\|_{\mathbb{H}^\lambda(\nLambda_t)}=\underset{\left(\theta_1,\ldots,\theta_n\right)\in
\nLambda_t}{\sup}\left|h\left(\theta_1,\ldots,\theta_n\right)\right|\\
&\quad \quad+\sum_{k=1}^n\sum_{i_1<\ldots <i_k=1}^n
\underset{0\leq\theta_{i_1}<
r_{i_1}\leq\theta_{i_2}<r_{i_2}\leq\ldots\leq\theta_{i_k}<r_{i_k}\leq
t}{\sup_{\left(\theta_1,\ldots,\theta_n\right)\in
\nLambda_t}}\frac{\left|\Delta^{i_1,\ldots,i_k}h\left(\theta_1,\ldots,\theta_n;r_{i_1},\ldots,r_{i_k}\right)\right|}
{\overset{k}{\underset{j=1}{\prod}}\left|r_{i_j}-\theta_{i_j}\right|^{\lambda}},
\end{align*}
where for $n=1$, $\Delta^1h(\theta,r)= h(r)-h(\theta)$, and for $n\ge 2$,
\begin{align*}
\Delta^{i_1}h\left(\theta_1,\ldots,\theta_n;r_{i_1}\right) &=
h\left(\theta_1,\ldots,r_{i_1},\ldots,\theta_n\right)
-h\left(\theta_1,\ldots,\theta_{i_1},\ldots,\theta_n\right),\\
\Delta^{i_1,\ldots,i_k}h\left(\theta_1,\ldots,\theta_n;r_{i_1},\ldots,r_{i_k}\right)&=
\Delta^{i_1,\ldots,i_{k-1}}h\left(\theta_1,\ldots,r_{i_k},\ldots,\theta_n;r_{i_1},\ldots,r_{i_{k-1}}\right)\\
&
-\Delta^{i_1,\ldots,i_{k-1}}h\left(\theta_1,\ldots,\theta_{i_k},\ldots,\theta_n;r_{i_1},\ldots,r_{i_{k-1}}\right),
\end{align*}
$2\leq k \leq n$.

Given a Banach space $\mathbb{B}$, we shall denote by $\mathcal{C}^\lambda\left([0,T]; \mathbb{B}\right)$
($\mathcal{C}^\lambda\left([0,T]\right)$, when $\mathbb{B}=\mathbb{R}$) the space of $\lambda$--H\"older continuous
functions endowed with the norm
\begin{equation*}
\Vert h\Vert_{\mathcal{C}^\lambda\left([0,T]; \mathbb{B}\right)}= \sup_{0\le t\le T}
\left\|h(t)\right\|_{\mathbb{B}}+\underset{s\neq t}{\sup_{0\leq s, t \leq T}}\frac{\left \Vert h(t)-h(s)\right\Vert_{\mathbb{B}}}{\left|t-s\right|^\beta}.
\end{equation*}
For $n=1$, the space $\mathbb{H}^\lambda(\Lambda_T^{(1)})$ is usually denoted by $\mathcal{C}^\lambda([0,T])$.
\medskip

For any $p\in[1,\infty[$, we denote by
$\mathbb{D}^{1,p}(\mathbb{B})$ the space of $\mathbb{B}$-valued
random variables $Y$ satisfying
\begin{equation*}
\left\|Y\right\|^p_{1,p,\mathbb{B}}:=\mathbb{E}\left(\left\|Y\right\|^p_\mathbb{B}\right)+\mathbb{E}\left(\int_{0}^{T}\left\|D_{\eta}Y\right\|^2_\mathbb{B}d\eta\right)^{\frac{p}{2}}<+\infty.
\end{equation*}
In the next Propositions $2.1$ and $2.2$, we shall use this
notation for $\mathbb{B}=\left|\mathcal{H}^H\right|$ and
$\mathbb{B}=\mathcal{H}_{K_H}$, respectively.

For $\mathbb{B}=\mathbb{R}$, we write $\mathbb{D}^{1,p}$ instead
of $\mathbb{D}^{1,p}(\mathbb{R})$, and this is the usual
Sobolev-Watanabe space associated with the Wiener process $W$.

Set
\begin{equation*}
X_t=\int_{0}^{t}u_\theta d\bh_\theta, \quad t\in \IntervalT,
\end{equation*}
where the meaning attached to the stochastic integral is that of \cite{AMN01}, that is
\begin{equation*}
\int_{0}^{t}u_\theta d\bh_\theta:= \int_0^t \left(K^\ast_{H,t}u\right) (\theta) dW_\theta.
\end{equation*}
The sample path properties of the multiple stochastic integrals investigated in this paper rely on
results telling us what properties on the integrand $u$ imply that the stochastic process $\{X_t, t\in[0,T]\}$ exist and its sample paths are H\"older
continuous functions.

An answer is provided by Propositions 3 and 1 of \cite{AMN01}, by
considering the fractional Brownian motion as Gaussian process and
taking in these statements $\alpha=H-\frac{1}{2}$, with $H\in
]\frac{1}{2}, 1[$, and $\alpha=\frac{1}{2}-H$, with
$H\in]0,\frac{1}{2}[$, respectively. The next Proposition
\ref{p1.1}, gives a slightly different result when $H\in
]\frac{1}{2}, 1[$, while Proposition \ref{p1.2} is just a
quotation of Proposition 1 of \cite{AMN01}.

\begin{prop}
\label{p1.1}
 Let $H\in]\frac{1}{2}, 1[$ and $p\in[2,\infty[$. Consider a stochastic process
$u=\left\{u_t,\;t \in \IntervalT\right\}$ belonging to
$\RegularSpaceIntegrand$. Then $u\in\SobolevSpace^{1,p}\left(\RegularSpace\right)$
and the stochastic integral process $X=\left\{X_t,  t \in \IntervalT\right\}$ is well
defined and consists of random variables in $L^p(\Omega)$.

Furthermore,
if $qH>1$, and $u \in L^{q}\left(\IntervalT;\SobolevSpace^{1,p}\right)$, then
$$\mathbb{E}|X_t-X_s|^p\le C|t-s|^{p(H-1/q)}.$$
Hence, if $p\left(H-1/q\right)>1,$  the process $X$ has
$\gamma$--H\"older continuous paths, a.s., with $\gamma \in ]0,H-1/q-1/p[$.
\end{prop}
{\it Proof}: By applying twice Minkowski's inequality, we have
\begin{align*}
&\mathbb{E}\left(\|u\|^p_{\RegularSpace}\right)\\
&=\left\Vert \int_0^T d\theta \left(\int_{\theta}^T \vert u_r\vert K(dr,\theta)\right)^2\right\Vert_{L^{\frac{p}{2}}(\Omega)}^{\frac{p}{2}}
\le \left(\int_0^T d\theta\left\Vert \int_{\theta}^T \vert u_r\vert K(dr,\theta)\right\Vert_{L^p(\Omega)}^2\right)^{\frac{p}{2}}\\
&\le \left(\int_0^T d\theta\left(\int_{\theta}^T dr \left\Vert u_r\Partialr\left(r,\theta\right)\right\Vert_{L^p(\Omega)}\right)^2\right)^{\frac{p}{2}}
\leq\left(\int_0^Td\theta \left(\int_\theta^T
\left\|u_r\right\|_{1,p}
K(dr,\theta)\right)^2\right)^{\frac{p}{2}}.
\end{align*}
By using first Fubini's theorem and  then Minkowski's inequality three times, we obtain
\begin{align*}
\mathbb{E}\left(\int_{0}^T\left\|D_\eta
u\right\|^2_{\RegularSpace}d\eta\right)^{\frac{p}{2}}\notag
&=\mathbb{E}\left(\int_{0}^Td\theta \left\Vert \int_\theta^T \vert D_{\cdot}u_r\vert K(dr,\theta)\right\Vert_{L^2([0,T])}^2\right)^{\frac{p}{2}}\\
&\le \left\Vert \int_{0}^Td\theta \left( \int_\theta^T \Vert D_{\cdot}u_r\Vert _{L^2([0,T])}K(dr,\theta)\right)^2\right\Vert_{L^{\frac{p}{2}}(\Omega)}^{\frac{p}{2}}\\
&\le \left(\int_{0}^Td\theta \left\Vert\int_\theta^T \Vert D_{\cdot}u_r\Vert _{L^2([0,T])}K(dr,\theta)\right\Vert_{L^p(\Omega)}^2\right)^{\frac{p}{2}}\\
&\le \left(\int_{0}^Td\theta \left(\int_\theta^T \Vert
u_r\Vert_{1,p}K(dr,\theta)\right)^2\right)^{\frac{p}{2}}.
\end{align*}
Since by (\ref{1.8})
\begin{equation*}
\int_0^T
\left(\int_\theta^T\left\|u_r\right\|_{1,p} K_H(dr,\theta)\right)^2 d\theta
\leq b_H^2 \left(\int_0^T
\left\|u_r\right\|_{1,p}^\frac{1}{H}dr\right)^{2H},
\end{equation*}
we obtain that $u\in \SobolevSpace^{1,p}\left(\RegularSpace\right)$.

For $0\le s<t$, following the proof of Proposition $1$ in \cite{AMN01}, we can write
\begin{equation*}
X_t-X_s=\int_0^s\left(\int_s^t u_r K_H
\left(dr,\theta\right)\right)dW_\theta+\int_s^t\left(\int_\theta^t
u_r K_H
\left(dr,\theta\right)\right)dW_\theta.
\end{equation*}
Thus, Meyer's inequality implies
\begin{equation*}
\mathbb{E}\left|X_t-X_s\right|^p\le  C(p)\left(S_1^2 +S_2^2\right)^{\frac{p}{2}},
\end{equation*}
with
\begin{align*}
S_1&= \left\|\, \Ind_{\left]0,s\right[}\int_ s^t \Vert
u_r\Vert_{1,p}\,\kh\left(dr,\cdot\right)\right\|_{L^2\left(\Intervalt\right)},\\
S_2&= \left\| \,\Ind_{\left]s,t\right[}\, \int_{\cdot}^t  \Vert
u_r\Vert_{1,p}\,\kh\left(dr,\cdot\right)\right\|_{L^2\left(\Intervalt\right)}.
\end{align*}
By comparing the domains of integration in the terms $S_1$ and $S_2$ above with that of
$\left\|\Vert u_r\Vert_{1,p} \;\Ind_{\left]s,t\right[}\right\|_{\left|\rkhsh\right|}$
we see that
\begin{equation*}
S_1^2 + S_2^2 \le \left\|\Vert u_r\Vert_{1,p}\;\,\Ind_{\left]s,t\right[}\right\|_{\left|\rkhsh\right|}^2.
\end{equation*}
Then, applying (\ref{1.8}) we obtain
\begin{equation*}
\mathbb{E}\left|X_t-X_s\right|^p
\leq C(p,H) \left\Vert \Vert u_r\Vert_{1,p}\;\, \Ind_{\left]s,t\right[}\right\Vert_{L^{\frac{1}{H}}([0,T])}^p.
\end{equation*}
Notice that this last integral is finite.

Fix $q\ge 1$ such that $qH>1$. By applying  H\"older's inequality
with $p'=qH$, $q'=\frac{qH}{qH-1}$, we reach
\begin{align*}
\left\Vert \Vert u_r\Vert_{1,p}\;\,
\Ind_{\left]s,t\right[}\right\Vert_{L^{\frac{1}{H}}([0,T])}^p&=
\left(\int_s^t
\left\|u_r\right\|_{1,p}^\frac{1}{H}dr\right)^{pH}\\
&\leq
\left(t-s\right)^{p\left(H-1/q\right)}\left(\int_s^t
\left\|u_r\right\|_{1,p}^q dr\right)^{\frac{p}{q}}\\
&\leq \left(t-s\right)^{p\left(H-1/q\right)}
\left\|u\right\|_{L^{q}
\left(\IntervalT;\SobolevSpace^{1,p}\right)}^p.
\end{align*}
Consequently,
$$
\mathbb{E}\left|X_t-X_s\right|^p\le C|t-s|^{p\left(H-1/q\right)},
$$
and we conclude by applying Kolmogorov's continuity criterion.
\hfill \qed

\begin{prop}
\label{p1.2} Let $H\in]0,\frac{1}{2}[$, $p\in[2,\infty[$. Suppose
that the stochastic process $u=\left\{u_t,\;t \in
\IntervalT\right\}$ belongs to $\SingularSpaceIntegrand$ for some
 $\lambda+H>\frac{1}{2}$. Then $u$ belongs to the
space
$\SobolevSpace^{1,p}\left(\SingularSpace\right)$ and
the stochastic integral process $\{X_t, t \in\IntervalT\}$ is well
defined; it consists of $L^p(\Omega)$ random variables and satisfies $\mathbb{E}|X_t-X_s|^p\le C|t-s|^{pH}$.
Consequently, if $pH>1$, a.s., the sample paths are
$\gamma$--H\"older continuous with $\gamma \in]0, H-1/p[$.
\end{prop}
\medskip


\section{Multiple stochastic integrals of deterministic functions with respect to a fractional Brownian motion}
\label{s2}
In this section, we study conditions on deterministic functions $h$ defined on $[0,T]^n$
allowing to define the indefinite multiple stochastic integral with respect to a fractional Brownian motion
\begin{equation*}
\fnint_{t}
\left(h\right)=\int_{\nLambda_t}h\left(\theta_1,\ldots,\theta_n\right)d\bh_{\theta_1}\cdots
d\bh_{\theta_n}, \, t\in[0,T].
\end{equation*}

{\bf $\bullet$ The case $H\in]\frac{1}{2},1[$.}
\medskip

Fix $0\leq s < t \leq T$. We set
\begin{align}
\left(K_{H,s,t}^{\ast,(1)}h\right)\left(\theta_1\right)&=\int_{\theta_1}^t
h\left(r_1\right) \kh\left(dr_1,\theta_1\right)\;
\Ind_{\left]s,t\right[}
\left(\theta_1\right)\notag\\
&\quad+ \int_ s^t h\left(r_1\right)
\kh\left(dr_1,\theta_1\right)\;\Ind_{\left]0,s\right[
}\left(\theta_1\right)\Ind_{\left\{s\neq0\right\}},\label{2.0}
\end{align}
and for any integer $n\ge 2$,
\begin{align}
&\left(\khst^{\ast,(n)}h\right)\left(\theta_1, \ldots,
\theta_n\right) \notag\\
&=\int_{\vee_{i=1}^n\theta_i}^t
\left(K_{H,r_n}^{\ast,(n-1)}h\left(\cdot,r_n\right)
\right)\left(\theta_1,\ldots,\theta_{n-1}\right) \kh\left(dr_n,\theta_n\right)\notag\\
&\quad\quad\times\Ind_{\left[0,t\right]^{n-1}}(\theta_1,\ldots,\theta_{n-1})\,\Ind_{]s,t[}(\theta_n)
\notag\\
&\quad+ \int_{\vee_{i=1}^{n-1}\theta_i\vee
s}^t\left(K_{H,r_n}^{\ast,(n-1)}h\left(\cdot,r_n\right)\right)\left(\theta_1,\ldots,\theta_{n-1}\right)
\kh\left(dr_n,\theta_n\right)\notag\\
&\quad\quad \times
\Ind_{\left[0,t\right]^{n-1}}(\theta_1,\ldots,\theta_{n-1})\,\Ind_{]0,s[}(\theta_n)\,\Ind_{\left\{s\neq0\right\}},\label{2.1}
\end{align}
where we write $K_{H,\cdot}^{\ast,(m)}$ for
$K_{H,0,\cdot}^{\ast,(m)}$.

\begin{prop}
\label{p2.1} Fix $H\in]\frac{1}{2}, 1[$ and a natural number $n\ge1$.
\begin{description}
\item{(a)} Let $h\in\nRegularSpace$. Then, for $0\leq s<t\leq T$,
\begin{equation}
\left\|\khst^{\ast,(n)}h\right\|_{L^2\left(\Intervalt^n\right)}\leq
2^{n/2}b_H^n\,
\left\|h\;\Ind_{\left]s,t\right[}\right\|_{L^\frac{1}{H}(\nLambda_t)},\label{2.1.0}
\end{equation}
where $b^H$ is the same constant as in
$(\ref{1.8})$. Thus,
$\khst^{\ast,(n)}:L^\frac{1}{H}(\nLambda_t)\longrightarrow
L^2\left(\Intervalt^n\right)$ defined in $(\ref{2.0})$ and
$(\ref{2.1})$ is a linear continuous operator.
\item{(b)}
If $h\in\nRSHolder$ with
$qH>1$, then for $0\leq s<t\leq T$,
\begin{equation}
\left\|\khst^{\ast,(n)}h\right\|_{L^2\left(\Intervalt^n\right)}\leq
2^{n/2}b_H^n\left\|h\right\|_{L^q(\nLambda_t)}\left|t-s\right|^{H-1/q}.\label{2.1.0.1}
\end{equation}
\end{description}
\end{prop}

{\it Proof}: Let $n=1$. Owing to $(\ref{2.0})$,
\begin{equation*}
\left\|\khst^{\ast,(1)}h\right\|_{L^2\left(\Intervalt\right)}\le T_1+T_2,
\end{equation*}
with
\begin{align*}
T_1&=
\left\| \,\Ind_{\left]s,t\right[}\, \int_{\cdot}^t h\left(r_1\right)\kh\left(dr_1,\cdot\right)\right\|_{L^2\left(\Intervalt\right)},\\
T_2&= \left\|\,
\Ind_{\left]0,s\right[}\;\Ind_{\left\{s\neq0\right\}}\int_ s^t
h\left(r_1\right)
\kh\left(dr_1,\cdot\right)\right\|_{L^2\left(\Intervalt\right)}.
\end{align*}
As in the proof of Proposition \ref{p1.1}, we see that
\begin{equation}
\label{clau}
T_1^2 + T_2^2 \le \left\|h\;\Ind_{\left]s,t\right[}\right\|_{\left|\rkhsh\right|}^2.
\end{equation}
Indeed, it suffices to compare the domains of integration of the terms $T_1$ and $T_2$ above with that of
$\left\|h\;\Ind_{\left]s,t\right[}\right\|_{\left|\rkhsh\right|}$.
Then, applying (\ref{1.8}) we obtain (\ref{2.1.0}) for $n=1$.

H\"older's
inequality with $p'=qH$, $q'=\frac{qH}{qH-1}$, yields
\begin{align*}
\left\|h\;\Ind_{\left]s,t\right[}\right\|_{L^{\frac{1}{H}}(\IntervalT)}\leq
\left\|h\;\Ind_{\left]s,t\right[}\right\|_{L^q(\IntervalT)}\left|t-s\right|^{H-1/q}\leq
\left\|h\right\|_{L^q(\IntervalT)}\left|t-s\right|^{H-1/q}.
\end{align*}
Hence $(\ref{2.1.0.1})$ holds for $n=1$.
\medskip

Assume now that $(\ref{2.1.0})$ holds up to an integer $n'\geq 1$.
Let $h \in L^\frac{1}{H}(\Lambda_t^{(n'+1)})$. Then
$h(\cdot,r_{n'+1})\in L^\frac{1}{H}(\Lambda^{(n')}_{r_{n'+1}})$,
for all $r_{n'+1}\in [0,t]$, a.e. and the induction hypothesis
yields
\begin{equation*}
\left\|K_{H,r_{n'+1}}^{\ast,(n')}h(\cdot,r_{n'+1})\right\|_{L^2\left(\left[0,r_{n'+1}\right]^{n'}\right)}\leq
2^{n'/2}b_H^{n'}\left\|h(\cdot,r_{n'+1})\right\|_{L^\frac{1}{H}(\Lambda^{(n')}_{r_{n'+1}})},
\end{equation*}
for all $r_{n'+1}\in [0,t]$, a.e.

From $(\ref{2.1})$ it follows that
$\left\|\khst^{\ast,(n'+1)}h\right\|^2_{L^2\left(\Intervalt^{n'+1}\right)}\leq
2\left(Q_1+Q_2\right)$, with

\begin{align*}
Q_1&=\int_s^t
d\theta_{n^{'}+1}\left\Vert\int_{\vee_{i=1}^{n'+1}\theta_i}^t
\left(K_{H,r_{n'+1}}^{\ast,(n')}h\left(\cdot,
r_{n'+1}\right)\right)(\cdot)
K_H(dr_{n'+1},\theta_{n'+1})\right\Vert_{L^2\left([0,t]^{n'}\right)}^2,\\
Q_2&=\int_0^s
d\theta_{n^{'}+1}\left\Vert\int_{\vee_{i=1}^{n'}\theta_i\vee s}^t
\left(K_{H,r_{n'+1}}^{\ast,(n')}h\left(\cdot,
r_{n'+1}\right)\right)(\cdot)
K_H(dr_{n'+1},\theta_{n'+1})\right\Vert_{L^2\left([0,t]^{n'}\right)}^2.
\end{align*}

Minkowski's inequality yields
\begin{align*}
Q_1&\leq \int_s^t d\theta_{n'+1}\left(\int_{\theta_{n'+1}}^t
\left\Vert K_{H,r_{n'+1}}^{\ast,(n')}h(\cdot,
r_{n'+1})\right\Vert_{L^2([0,r_{n'+1}]^{n'})} K_H(dr_{n'+1},\theta_{n'+1})\right)^2\\
&\leq  2^{n'} b_H^{2n'} \int_s^t  d\theta_{n'+1}\left(\int_{\theta_{n'+1}}^t
\left\|h(\cdot,r_{n'+1})\right\|_{L^\frac{1}{H}(\Lambda^{(n')}_{r_{n'+1}})}
K_H\left(dr_{n'+1},\theta_{n'+1}\right)\right)^2
\end{align*}
and
\begin{align*}
Q_2&\leq \int_0^s d\theta_{n'+1}\left(\int_s^t \left\Vert
K_{H,r_{n'+1}}^{\ast,(n')}h(\cdot,
r_{n'+1})\right\Vert_{L^2([0,r_{n'+1}]^{n'})} K_H(dr_{n'+1},\theta_{n'+1})\right)^2\\
&\leq 2^{n'} b_H^{2n'}\int_0^s  d\theta_{n'+1}\left(\int_s^t
\left\|h(\cdot,r_{n'+1})\right\|_{L^\frac{1}{H}(\Lambda^{(n')}_{r_{n'+1}})}
K_H\left(dr_{n'+1},\theta_{n'+1}\right)\right)^2.
\end{align*}
Notice that, except for the constant $2^{n'} b_H^{2n'}$, the upper
bounds of the terms $Q_1$, $Q_2$ coincide with $T_1^2$ and
$T_2^2$, respectively, with $h:=
\left\|h(\cdot,r_{n'+1})\right\|_{L^\frac{1}{H}(\Lambda^{(n')}_{r_{n'+1}})}$.

Thus, using (\ref{clau}) and (\ref{1.8})
\begin{align*}
\left\|\khst^{\ast,(n'+1)} h\right\|_{L^2\left(\Intervalt^{n'+1}\right)}&\le 2^{\frac{1}{2}} \left(Q_1+Q_2\right)^{\frac{1}{2}}\\
&\le 2^{\frac{n'+1}{2}} b_H^{n'}\left\|\, \left\|h(\cdot,r_{n'+1})\right\|_{L^{\frac{1}{H}}(\Lambda^{(n')}_{r_{n'+1}})}
   \Ind_{\left]s,t\right[}\right\|_{\vert\mathcal{H}^H\vert}\\
&\le  2^{\frac{n'+1}{2}} b_H^{n'+1} \left\|\, \left\|h(\cdot,r_{n'+1})\right\|_{L^{\frac{1}{H}}(\Lambda^{(n')}_{r_{n'+1}})}
   \Ind_{\left]s,t\right[}\right\|_{L^{\frac{1}{H}}([0,T])}\\
   &\le 2^{\frac{n'+1}{2}} b_H^{n'+1} \left\Vert h\, \Ind_{\left]s,t\right[}\right\Vert_{L^{\frac{1}{H}}(\Lambda_t^{(n'+1)})}.
\end{align*}

This proves (\ref{2.1.0}) for $n'+1$.

The upper bound (\ref{2.1.0.1}) for $n'+1$ follows by applying
H\"older's inequality, as we did for $n=1$ in the first step of
the proof. \hfill \qed

\begin{thm}\quad
\label{t2.1}
\begin{description}
\item{(A)} With the same hypothesis as in $(a)$ of Proposition
$\ref{p2.1}$, the integral stochastic process $\fnint(h)=
\left\{\fnint_t \left(h\right);\;t \in \IntervalT \right\}$, $n\ge 1$,  is
well defined as an iterated integral and for any $0\leq s<t \leq T$,
\begin{equation}
\fnint_t \left(h\right)-\fnint_s
\left(h\right)=\int_{\Intervalt^n}\left(\khst^{\ast,(n)}h\right)\left(\theta_1,\ldots,\theta_n\right)dW_{\theta_1}\cdots
dW_{\theta_n},\label{2.1.2}
\end{equation}
with $\khst^{\ast,(n)}$, given in $(\ref{2.0})$ and $(\ref{2.1})$.

\item{(B)} Suppose the same hypotheses as in $(b)$ of Proposition $\ref{p2.1}$.
Then, for any $p\in [2,\infty[$ and $0\le s<t\le T$,
\begin{align}
\left\|\fnint_t \left(h\right)-\fnint_s
\left(h\right)\right\|_{L^p(\Omega)}\leq
C\left\Vert K^{\ast,(n)}_{H,s,t}\right\Vert_{L^2([0,t]^n)}\le
C\left|t-s\right|^{H-1/q},\label{2.1.3}
\end{align}
for some positive constant $C$ depending on $q$, $p$, $h$ and $H$.

Consequently, the sample paths of $\fnint \left(h\right)$ are
$\gamma$--H\"older continuous with $\gamma \in
\left]0,H-1/q\right[$.
\end{description}
\end{thm}

{\it Proof}: Let us show $(A)$. We start by
noticing that for $n=1$ the equality $(\ref{2.1.2})$ has already been
met in Proposition $\ref{p1.1}$ (see the proof of
Proposition $3$ in \cite{AMN01}).

Assume that
$(\ref{2.1.2})$ holds true up to an integer $n'\geq 1$. Let
$h\in L^{\frac{1}{H}}(\Lambda_T^{(n'+1)})$. By the induction
assumption, for any $r_{n'+1} \in \Intervalt$, a.e., the random
variable
$\MultipleIntegral^{(n'),H}_{r_{n'+1}}\left(h\left(\cdot,r_{n'+1}\right)\right)$
is well defined as an iterated integral.

Fix $p \geq 2$. The rules of the Malliavin derivative for the standard Brownian motion and the hypercontractivity
inequality (see for instance \cite{LT91}) imply
\begin{align*}
\left\|\MultipleIntegral^{(n'),H}_{r_{n'+1}}\left(h\left(\cdot,r_{n'+1}\right)\right)\right\|_{1,p}
&\leq
C(n',p)\left\|K_{H,r_{n'+1}}^{\ast,(n')}h(\cdot,r_{n'+1})\right\|_{L^2\left(\left[0,r_{n'+1}\right]^{n'}\right)}\\
&\leq
C(n',p,H)\left\|h(\cdot,r_{n'+1})\right\|_{L^\frac{1}{H}(\Lambda^{(n')}_{r_{n'+1}})},
\end{align*}
for any $r_{n'+1}\in \Intervalt$, a.e.

Thus,
\begin{align*}
&\left(\int_0^T\left\|\MultipleIntegral^{(n'),H}_{r_{n'+1}}\left(h\left(\cdot,r_{n'+1}\right)\right)\right\|_{1,p}^\frac{1}{H}dr_{n'+1}\right)^H\notag\\
&\leq
C(n',p,H)\left(\int_0^T\left\|h(\cdot,r_{n'+1})\right\|_{L^\frac{1}{H}(\Lambda^{(n')}_{r_{n'+1}})}^\frac{1}{H}dr_{n'+1}\right)^H\notag\\
&=C(n',p,H)\left\|h\right\|_{L^\frac{1}{H}(\Lambda^{(n'+1)}_T)}<\infty.\notag
\end{align*}
Thus, the process $\left\{\MultipleIntegral^{(n'),H}_{r_{n'+1}}\left(h\left(\cdot,r_{n'+1}\right)\right);\;r_{n'+1}\in \IntervalT\right\}$
belongs to $L^\frac{1}{H}\left(\IntervalT;\;\SobolevSpace^{1,p}\right)$, for
any $p \geq 2$.

By applying Proposition \ref{p1.1}, we can write
\begin{equation*}
\MultipleIntegral^{(n'+1),H}_{t}\left(h\right)-\MultipleIntegral^{(n'+1),H}_{s}\left(h\right)
=\int_0^T \left(K_{H,s,t}^{\ast,(1)}\MultipleIntegral^{(n'),H}_{\star}\left(h(\cdot,\star)\right)\right)(\theta_{n'+1})dW_{\theta_{n'+1}}.
\end{equation*}
Owing to (\ref{2.0}), the last term can be decomposed into the sum
of two terms denoted by $M_1$, $M_2$, and defined as follows:
\begin{align*}
&M_1=\int_s^t\Big(\int_{\theta_{n'+1}}^t\MultipleIntegral^{(n'),H}_{r_{n'+1}}
\left(h(\cdot,r_{n'+1})\right)K_H(dr_{n'+1},\theta_{n'+1})\Big)dW_{\theta_{n'+1}},\\
&M_2=\int_0^s \Big(\int_s^t\MultipleIntegral^{(n'),H}_{r_{n'+1}}
\left(h(\cdot,r_{n'+1})\right)K_H(dr_{n'+1},\theta_{n'+1})\;\Ind_{\left\{s\neq0\right\}}\Big)dW_{\theta_{n'+1}}.
\end{align*}

We now apply the induction assumption to write all these terms by
means of $(n'+1)$--multiple integrals, obtaining
\begin{align}
&M_1=\int_s^t\Big(\int_{\theta_{n'+1}}^t\Big(\int_{\left[0,r_{n'+1}\right]^{n'}}
\left(K_{H,r_{n'+1}}^{\ast,(n')}
h\left(\cdot,r_{n'+1}\right)\right)\left(\theta_1,\ldots,\theta_{n'}\right)dW_{\theta_1}\cdots
dW_{\theta_{n'}}
\Big)\notag\\
&\quad \quad \times K_H(dr_{n'+1},\theta_{n'+1}) \Big)dW_{\theta_{n'+1}},\label{2.1.4}\\
&M_2=\int_0^s
\Big(\int_s^t\Big(\int_{\left[0,r_{n'+1}\right]^{n'}}\left(K_{H,r_{n'+1}}^{\ast,(n')}
h\left(\cdot,r_{n'+1}\right)\right)\left(\theta_1,\ldots,\theta_{n'}\right)dW_{\theta_1}\cdots
dW_{\theta_{n'}} \Big)\notag\\
&\quad \quad \times
K_H(dr_{n'+1},\theta_{n'+1})\;\Ind_{\left\{s\neq0\right\}}\Big)dW_{\theta_{n'+1}}.\label{2.1.5}
\end{align}
Finally, by the stochastic Fubini theorem applied to
$(\ref{2.1.4})$ and $(\ref{2.1.5})$, we obtain
\begin{align*}
&\MultipleIntegral^{(n'+1),H}_{t}(h) - \MultipleIntegral^{(n'+1),H}_{s}(h)\\
&=\int_{[0,t]^{n'}\times]s,t[}\Big(\int_{\vee_{i=1}^{n'+1}\theta_i}^t
\left(K_{H,r_{n'+1}}^{\ast,(n')}
h\left(\cdot,r_{n'+1}\right)\right)\left(\theta_1,\dots,\theta_{n'}\right)\\
&\quad\times
K_H(dr_{n'+1},\theta_{n'+1})\Big)dW_{\theta_1}\cdots
dW_{\theta_{n'}}dW_{\theta_{n'+1}}\\
&\quad+\int_{\left[0,t\right]^{n'}\times\left]0,s\right[}
\Big(\int_{\vee_{i=1}^{n'}\theta_i\vee
s}^t\left(K_{H,r_{n'+1}}^{\ast,(n')}
h\left(\cdot,r_{n'+1}\right)\right)\left(\theta_1,\ldots,\theta_{n'}\right)\\
&\quad \quad \times
K_H(dr_{n'+1},\theta_{n'+1})\;\Ind_{\left\{s\neq0\right\}}\Big)dW_{\theta_1}\cdots
dW_{\theta_{n'}}dW_{\theta_{n'+1}}.
\end{align*}
That is,
\begin{equation*}
\MultipleIntegral^{(n'+1),H}_{t}
\left(h\right)-\MultipleIntegral^{(n'+1),H}_{s} \left(h\right)
=\int_{\Intervalt^{n'+1}}\left(\khst^{\ast,(n'+1)}h\right)
\left(\theta_1,\ldots,\theta_{n'+1}\right)dW_{\theta_1}\cdots
dW_{\theta_{n'+1}}
\end{equation*}
(see (\ref{2.1})). This ends the proof of (\ref{2.1.2}) and also that of $(A)$.

For the proof of $(B)$, we consider (\ref{2.1.2}) and we apply
the hypercontractivity inequality and (\ref{2.1.0.1}). Then the
conclusion on sample path regularity is a consequence of the
Kolmogorov's criterion. \hfill\qed

\medskip

 In the next section, we shall consider indefinite multiple integrals with integrands depending on the upper bound of the
 integration domain for which we shall apply the following Corollary.
 \begin{cor}
\label{c2.1} Let $H\in]\frac{1}{2},1[$, $\beta\in ]0,1[$ and $h$
be a measurable function defined on $[0,T]^{n+1}$ such that
 the mapping $t\longmapsto h(\cdot,t)$
belongs to
$\mathcal{C}^\beta\left([0,T];L^q(\Lambda_T^{(n)})\right)$, for
$qH>1$. Then, the integral process
$\left\{\fnint_t\left(h(\cdot,t)\right),\;t \in\IntervalT\right\}$
given in Theorem $\ref{t2.1}$ has $\gamma$--H\"older continuous
sample paths, a.s., with $\gamma\in\left]0,\beta \wedge
\left(H-1/q\right)\right[$.
\end{cor}
{\it Proof}: Let $p \in [2,\infty[$. The hypercontractivity inequality and
(\ref{2.1.0.1}) yield
\begin{align*}
&\left\|\MultipleIntegral^{(n),H}_{t}
\left(h(\cdot,t)\right)-\MultipleIntegral^{(n),H}_{s}
\left(h(\cdot,s)\right)\right\|_{L^p(\Omega)}\notag\\
&\leq
\left\|\int_{\Intervalt^n}\left(\khst^{\ast,(n)}h\left(\cdot,t\right)\right)
\left(\theta_1,\ldots,\theta_n\right)dW_{\theta_1}\cdots
dW_{\theta_n}\right\|_{L^p(\Omega)}\notag\\
&\quad+\left\|\int_{\left[0,s\right]^n}\left(K_{H,s}^{\ast,(n)}\left[\Delta^{n+1}h\left(\cdot,s;t\right)\right]\right)
\left(\theta_1,\ldots,\theta_n\right)dW_{\theta_1}\cdots
dW_{\theta_n}\right\|_{L^p(\Omega)}\notag\\
&\leq
C(n,p)\left[\left\|\khst^{\ast,(n)}h\left(\cdot,t\right)\right\|_{L^{2}\left(\Intervalt^{n}\right)}+\left\|K_{H,s}^{\ast,(n)}
\left[\Delta^{n+1}h\left(\cdot,s;t\right)\right]\right\|_{L^2\left(\left[0,s\right]^{n}\right)}\right]\notag\\
&\leq
C(n,p,T,H)\left[\left\|h\right\|_{\mathcal{C}^\beta\left([0,T];L^q(\Lambda_T^{(n)})\right)}\left|t-s\right|^{H-1/q}
+\left\|h\right\|_{\mathcal{C}^\beta\left([0,T];L^q(\Lambda_T^{(n)})\right)}\left|t-s\right|^\beta\right] \notag\\
&\leq
C(n,p,T,H)\left\|h\right\|_{\mathcal{C}^\beta\left([0,T];L^q(\Lambda_T^{(n)})\right)}\left|t-s\right|^{\left(H-1/q\right)
\wedge \beta}.\notag
\end{align*}
The conclusion follows by applying Kolmogorov's criterion.
\hfill\qed

\medskip

{\bf $\bullet$ The case  $H\in]\frac{1}{4},\frac{1}{2}[$.}
\medskip

Let us introduce the functions that will appear as kernels of
increments of the indefinite multiple integrals. Fix $0\le s<t\le
T$. We set
\begin{align}
\left(\khst^{\ast,(1)}h\right)\left(\theta_1\right) \notag
&=
h\left(\theta_1\right)\kh\left(t,\theta_1\right)
\;\Ind_{\left]s,t\right[}\left(\theta_1\right)\\
&\quad+\int_{\theta_1}^t
\left[h\left(r_1\right)-h\left(\theta_1\right)\right]
\kh\left(dr_1,\theta_1\right)\; \Ind_{\left]s,t\right[}
\left(\theta_1\right)\notag\\
&\quad+
\int_ s^th\left(r_1\right)
\kh\left(dr_1,\theta_1\right)\;\Ind_{\left]0,s\right[
}\left(\theta_1\right)\,\Ind _{\left\{s\neq0\right\}},\label{u}
\end{align}
and for any integer $n\ge 2$,
\begin{align}
&\left(\khst^{\ast,(n)}h\right)\left(\theta_1, \ldots,
\theta_n\right) \notag\\
&=
\left(K_{H,\theta_n}^{\ast,(n-1)}h\left(\cdot,\theta_n\right)\right)\left(\theta_1,\ldots,\theta_{n-1}\right)\kh\left(t,\theta_n\right)
\;\Ind_{\left[0,\theta_n\right]^{n-1}}(\theta_1,\ldots,\theta_{n-1})\,\Ind_{]s,t[}(\theta_n)
\notag\\
&\quad
+\int_{\theta_n}^t
\left(K_{H,\theta_n}^{\ast,(n-1)}\left[h\left(\cdot,r_n\right)
-h\left(\cdot,\theta_n\right)\right]\right)\left(\theta_1,\ldots,\theta_{n-1}\right)
 \kh\left(dr_n,\theta_n\right)\notag\\
&\quad\quad\times\Ind_{\left[0,\theta_n\right]^{n-1}}(\theta_1,\ldots,\theta_{n-1})\,\Ind_{]s,t[}(\theta_n)
\notag\\
&\quad+
\int_{\vee_{i=1}^n\theta_i}^t\left(K_{H,\theta_n,r_n}^{\ast,(n-1)}h\left(\cdot,r_n\right)\right)\left(\theta_1,\ldots,\theta_{n-1}\right)
\kh\left(dr_n,\theta_n\right)\notag\\
&\quad \quad \times
\Ind_{\left[0,t\right]^{n-1}}(\theta_1,\ldots,\theta_{n-1})\,\Ind_{]s,t[}(\theta_n)\notag\\
&\quad+
\int_{\vee_{i=1}^{n-1}\theta_i\vee
s}^t\left(K_{H,r_n}^{\ast,(n-1)}h\left(\cdot,r_n\right)\right)\left(\theta_1,\ldots,\theta_{n-1}\right)
\kh\left(dr_n,\theta_n\right)\notag\\
&\quad\quad \times
\Ind_{\left[0,t\right]^{n-1}}(\theta_1,\ldots,\theta_{n-1})\,\Ind_{]0,s[}(\theta_n)\;\Ind_{\left\{s\neq0\right\}},\label{2.2.2}
\end{align}
where we write $K_{H,\cdot}^{\ast,(m)}$ instead of $K_{H,0,\cdot}^{\ast,(m)}$.

\begin{prop}
\label{p2.2} Fix $H\in]\frac{1}{4}, \frac{1}{2}[$ and a natural
number $n\ge 1$. Let $h\in\nSingularSpace$, for some $\lambda$
satisfying $\lambda+H>\frac{1}{2}$.  Then, for $0\leq s<t\leq T$,
\begin{equation}
\left\|\khst^{\ast,(n)}h\right\|_{L^2\left(\Intervalt^n\right)}\leq
C(T,\lambda,H)\left\|h\right\|_{\mathbb{H}^\lambda(\nLambda_t)}\left|t-s\right|^{H},\label{2.2.4}
\end{equation}
with some positive constant $C(T,\lambda,H)$. Thus,
$\khst^{\ast,(n)}:\mathbb{H}^\lambda(\nLambda_t)\longrightarrow
L^2\left(\Intervalt^n\right)$ defined in $(\ref{u})$ and
$(\ref{2.2.2})$ is a linear continuous operator.
\end{prop}

{\it Proof}: Let $n=1$. Owing to (\ref{u}),
\begin{align*}
&\left\|\khst^{\ast,(1)} h\right\|_{L^{2}\left(\Intervalt\right)}\notag\\
&\leq\left\|h\left(\cdot\right)\kh(t,\cdot)\;\Ind_{\left]s,t\right[}\right\|
_{L^{2}\left(\Intervalt\right)}+\left\|\int_{\cdot}^t
\left[\Delta^1h\left(\cdot;r_1\right)\right] K_H(dr_1,\cdot)\;\Ind_{\left]s,t\right[}\right\|_{L^{2}\left(\Intervalt\right)}\notag\\
&\quad+\left\|\int_s^t
h\left(r_1\right)K_H\left(dr_1,\cdot\right)\;\Ind_{\left]0,s\right[}\;\Ind
_{\left\{s\neq0\right\}}\right\|_{L^{2}\left(\Intervalt\right)}\notag\\
&\leq \left\|h\right\|_{\mathbb{H}^\lambda(\Lambda^{(1)}_t)}
\Big[\left\|\kh(t,\cdot)\;\Ind_{\left]s,t\right[}\right\|_{L^{2}\left(\Intervalt\right)}
+ \left\Vert\int_{\cdot}^{t}\left\vert
r_1-\cdot\right\vert^\lambda\left\vert
K_H(dr_1,\cdot)\right\vert\;\Ind_{\left]s,t\right[}\right\Vert
_{L^{2}\left(\Intervalt\right)}\notag\\
&\quad+ \left\|\int_{s}^{t}\left\vert
K_H\left(dr_1,\cdot\right)\right\vert\;\Ind_{\left]0,s\right[}\;\Ind
_{\left\{s\neq0\right\}}\right\|
_{L^{2}\left(\Intervalt\right)}\Big]\\
&\leq C(T,\lambda,H)
\left\|h\right\|_{\mathbb{H}^\lambda(\Lambda^{(1)}_t)}|t-s|^H,
\end{align*}
where in the last estimate we have used $(\ref{1.20})$. Thus,
$(\ref{2.2.4})$ holds for $n=1$.

\medskip
Suppose now  that $(\ref{2.2.4})$ holds up to an integer $n'\ge
1$. Let $h\in\mathbb{H}^{\lambda}(\Lambda^{(n'+1)}_t)$. The
functions $h(\cdot,\eta_{n'+1})$ and
$\Delta^{n'+1}h(\cdot,\eta_{n'+1};\tau_{n'+1})$, for any fixed
$\eta_{n'+1}\le t$ and $\tau_{n'+1}<\eta_{n'+1}\le t$,
respectively, belong to
$\mathbb{H}^{\lambda}(\Lambda^{(n')}_{\eta_{n'+1}})$.
From $(\ref{2.2.2})$ it follows that $\left\|\khst^{\ast,(n'+1)}h\right\|^2_{L^2\left(\Intervalt^{n'+1}\right)}\leq
2^3\left(\sum_{i=1}^4 R_i\right)$,
with

\begin{align*}
&R_1=
\int_{\left[0,\theta_{n'+1}\right]^{n'}\times\left]s,t\right[}\left(K_{H,\theta_{n'+1}}^{\ast,(n')}h(\cdot,\theta_{n'+1})\right)^2
\left(\theta_1,\ldots,\theta_{n'}\right)\\
&\quad \quad\times \kh\left(t,\theta_{n'+1}\right)^2d\theta_1
\cdots d\theta_{n'}d\theta_{n'+1};\notag\\
&R_2=\int_{\left[0,\theta_{n'+1}\right]^{n'}\times
\left]s,t\right[}\Big(\int_{\theta_{n'+1}}^t
\left(K_{H,\theta_{n'+1}}^{\ast,(n')}\left[\Delta^{n'+1}h
\left(\cdot,\theta_{n'+1};r_{n'+1}\right)\right]\right)\left(\theta_1,\ldots,\theta_{n'}\right)\\
&\quad \quad\times
\left|K_{H}\left(dr_{n'+1},\theta_{n'+1}\right)\right|\Big)^2d\theta_{1}\cdots
d\theta_{n'}d\theta_{n'+1};\\
&R_3=\int_{\left[0,t\right]^{n'}\times
\left]s,t\right[}\Big(\int_{\vee_{i=1}^{n'+1}\theta_i}^t
\left(K_{H,\theta_{n'+1},r_{n'+1}}^{\ast,(n')}
h\left(\cdot,r_{n'+1}\right)\right)\left(\theta_1,\ldots,\theta_{n'}\right)\\
&\quad\quad\times\left|K_{H}\left(dr_{n'+1},\theta_{n'+1}\right)\right|\Big)^2
d\theta_{1}\cdots
d\theta_{n'}d\theta_{n'+1};\\
&R_4=\int_{\left[0,t\right]^{n'}\times
\left]0,s\right[}\Big(\int_{\vee_{i=1}^{n'}\theta_i\vee s}^t
\left(K_{H,r_{n'+1}}^{\ast,(n')}h\left(\cdot,r_{n'+1}\right)\right)\left(\theta_1,\ldots,\theta_{n'}\right)\\
&\quad \quad\times
\left|K_{H}\left(dr_{n'+1},\theta_{n'+1}\right)\right|
\Big)^2 d\theta_{1}\cdots d\theta_{n'}d\theta_{n'+1}.
\end{align*}
Applying Minkowski's inequality to each one of these terms and
the induction assumption, we obtain,
\begin{align}
&R_1\leq\int_s^t \left\|K_{H,\theta_{n'+1}}^{\ast,(n')}h
\left(\cdot,\theta_{n'+1}\right)\right\|_{L^2\left(\left[0,\theta_{n'+1}\right]^{n'}\right)}^2
\kh\left(t,\theta_{n'+1}\right)^2d\theta_{n'+1}\nonumber\\
&\quad\leq C(T,\lambda,H) \int_s^t
 \left\|h(\cdot,\theta_{n'+1})\right\|^2_{\mathbb{H}^\lambda(\Lambda^{(n')}_{\theta_{n'+1}})}\kh\left(t,\theta_{n'+1}\right)^2d\theta_{n'+1}\nonumber\\
&\quad \leq C(T,\lambda,H)
\left\|h\right\|^2_{\mathbb{H}^\lambda(\Lambda_t^{(n'+1)})}|t-s|^{2H},\label{A}
\end{align}
\begin{align}
&R_2\leq\int_s^t
\Big(\int_{\theta_{n'+1}}^t\left\|K_{H,\theta_{n'+1}}^{\ast,(n')}
\left[\Delta^{n'+1}h
\left(\cdot,\theta_{n'+1};r_{n'+1}\right)\right]\right\|_{L^2\left(\left[0,\theta_{n'+1}\right]^{n'}\right)}\nonumber\\
&\quad\quad \times \left| K_{H}\left(dr_{n'+1},\theta_{n'+1}\right)\right|\Big)^2d\theta_{n'+1}\nonumber\\
&\quad\leq C(T,\lambda,H)
 \int_s^t\Big(\int_{\theta_{n'+1}}^t
 \left\|\Delta^{n'+1}h(\cdot,\theta_{n'+1};r_{n'+1})\right\|_{\mathbb{H}^\lambda(\Lambda^{(n')}_{\theta_{n'+1}})}\nonumber\\
&\quad \quad\times\left| K_{H}\left(dr_{n'+1},\theta_{n'+1}\right)\right|\Big)^2d\theta_{n'+1}\nonumber\\
&\quad\leq C(T,\lambda,H)\left\|h\right\|^2_{\mathbb{H}^\lambda(\Lambda_t^{(n'+1)})}\nonumber\\
&\quad \quad \times \int_s^t\Big(\int_{\theta_{n'+1}}^t
\left |r_{n'+1}-\theta_{n'+1}\right|^\lambda\left| K_{H}(dr_{n'+1},\theta_{n'+1})\right|\Big)^2d\theta_{n'+1}\nonumber\\
&\quad \leq
C(T,\lambda,H)\left\|h\right\|^2_{\mathbb{H}^\lambda(\Lambda_t^{(n'+1)})}\left|t-s\right|^{2(H+\lambda)},\label{B}
\\\nonumber
\\
&R_3\leq\int_s^t\Big(\int_{\theta_{n'+1}}^t\left\|K_{H,\theta_{n'+1},r_{n'+1}}^{\ast,(n')}
h\left(\cdot,r_{n'+1}\right)\right\|_{L^2\left(\left[0,r_{n'+1}\right]^{n'}\right)}\nonumber\\
&\quad \quad\times\left| K_{H}
\left(dr_{n'+1},\theta_{n'+1}\right)\right|\Big)^2
d\theta_{n'+1}\nonumber\\
&\quad \leq C(T,\lambda,H)\int_s^t\Big(\int_{\theta_{n'+1}}^t
\left\|h(\cdot,r_{n'+1})\right\|_{\mathbb{H}^\lambda(\Lambda^{(n')}_{r_{n'+1}})}\nonumber\\
&\quad \quad\times\left|r_{n'+1}-\theta_{n'+1}\right|^H
\left| K_{H}\left(dr_{n'+1},\theta_{n'+1}\right)\right|\Big)^2
d\theta_{n'+1}\nonumber\\
&\quad\leq C(T,\lambda,H)\left\|h\right\|^2_{\mathbb{H}^\lambda(\Lambda_t^{(n'+1)})}\nonumber\\
&\quad \quad \times\int_s^t\Big(\int_{\theta_{n'+1}}^t
\left|r_{n'+1}-\theta_{n'+1}\right|^H \left|
K_{H}(dr_{n'+1},\theta_{n'+1})\right|\Big)^2
d\theta_{n'+1}\nonumber\\
&\quad \leq
C(T,\lambda,H)\left\|h\right\|^2_{\mathbb{H}^\lambda(\Lambda_t^{(n'+1)})}\left|t-s\right|^{4H},\label{C}
\\\nonumber
\\
&R_4\leq\int_0^s\Big(\int_s^t\left\|K_{H,r_{n'+1}}^{\ast,(n')}
h\left(\cdot,r_{n'+1}\right)\right\|_{L^2\left(\left[0,r_{n'+1}\right]^{n'}\right)}\nonumber\\
&\quad \quad\times
\left|K_{H}\left(dr_{n'+1},\theta_{n'+1}\right)\right|
\Big)^2 d\theta_{n'+1}\nonumber\\
&\quad \leq C(T,\lambda,H)
\int_0^s\Big(\int_s^t\left\|h(\cdot,r_{n'+1})\right\|_{\mathbb{H}^\lambda(\Lambda^{(n')}_{r_{n'+1}})}\left|
K_{H}\left(dr_{n'+1},\theta_{n'+1}\right)\right|
\Big)^2d\theta_{n'+1}\nonumber\\
&\quad \leq
C(T,\lambda,H)\left\|h\right\|^2_{\mathbb{H}^\lambda(\Lambda_t^{(n'+1)})}\int_0^s\left(\int_s^t\left\vert
K_{H}(dr_{n'+1},\theta_{n'+1})\right\vert\right)^2d\theta_{n'+1}\nonumber\\
&\quad \leq
C(T,\lambda,H)\left\|h\right\|^2_{\mathbb{H}^\lambda(\Lambda_t^{(n'+1)})}\left|t-s\right|^{2H}.\label{D}
\end{align}
With (\ref{A}), (\ref{B}), (\ref{C}), (\ref{D}), we see that
(\ref{2.2.4}) holds for $n=n'+1$ and this ends the
proof of this proposition. \hfill \qed

\begin{thm}
\label{t2.2} With the same hypotheses as in Proposition $\ref{p2.2}$, the indefinite integral stochastic process $\fnint
\left(h\right)=\left\{\fnint_t \left(h\right), \;t \in \IntervalT
\right\}$ is well defined as an iterated integral. Moreover, for
any $0\le s<t\le T$,
\begin{equation}
\fnint_t \left(h\right)-\fnint_s
\left(h\right)=\int_{\Intervalt^n}\left(\khst^{\ast,(n)}h\right)
\left(\theta_1,\ldots,\theta_n\right)dW_{\theta_1}\cdots
dW_{\theta_n},\label{2.5}
\end{equation}
with $\khst^{\ast,(n)}$, given in $(\ref{u})$ and
$(\ref{2.2.2})$.

Thus, for any $p\in[2,\infty[$,
\begin{equation}
\label{2.5.0.0}
\left\Vert \fnint_t\left(h\right)-\fnint_s\left(h\right)\right\Vert_{L^p(\Omega)}\le
C \left\Vert K_{H,s,t}^{\ast,(n)}h\right\Vert_{L^2([0,t]^n)}\le
C |t-s|^H,
\end{equation}
for some positive constant $C$ depending on $p$, $h$, $T$, $\lambda$ and $H$.

Consequently, the sample paths of $\fnint \left(h\right)$ are
$\gamma$--H\"older continuous with $\gamma\in]0,H[$.
\end{thm}

{\it Proof}: Let us prove first $(\ref{2.5})$. For $n=1$, it is an
immediate consequence of Proposition \ref{p1.2}. The formula
$(\ref{2.5})$ is given in the proof of Proposition 1 of
\cite{AMN01}.

Assume that $(\ref{2.5})$ holds up to an integer $n'\ge 1$.
Consider a function $h\in\mathbb{H}^\lambda(\Lambda^{(n'+1)}_T)$.
For any $0\le s< t \le T$, we can write
\begin{align}
&\MultipleIntegral^{(n'),H}_{t}
\left(h\left(\cdot,t\right)\right)-\MultipleIntegral^{(n'),H}_{s}
\left(h\left(\cdot,s\right)\right)\nonumber\\
&\quad=\MultipleIntegral^{(n'),H}_{t}
\left(h\left(\cdot,t\right)\right)-\MultipleIntegral^{(n'),H}_{s}
\left(h\left(\cdot,t\right)\right)+\MultipleIntegral^{(n'),H}_{s}
\left(\left[\Delta^{n'+1}h\left(\cdot,s;t\right)\right]\right)\label{50}.
\end{align}
By the induction assumption, the following representations hold:
\begin{align*}
&\MultipleIntegral^{(n'),H}_{t}
\left(h\left(\cdot,t\right)\right)-\MultipleIntegral^{(n'),H}_{s}
\left(h\left(\cdot,t\right)\right)\\
&=\int_{\Intervalt^{n'}}\left(\khst^{\ast,(n')}h\left(\cdot,t\right)\right)
\left(\theta_1,\ldots,\theta_{n'}\right)dW_{\theta_1}\cdots
dW_{\theta_{n'}},\\
&\MultipleIntegral^{(n'),H}_{s}
\left(\left[\Delta^{n'+1}h\left(\cdot,s;t\right)\right]\right)\\
&=\int_{\left[0,s\right]^{n'}}\left(K_{H,s}^{\ast,(n')}\left[\Delta^{n'+1}h\left(\cdot,s;t\right)\right]\right)
\left(\theta_1,\ldots,\theta_{n'}\right)dW_{\theta_1}\cdots
dW_{\theta_{n'}}.
\end{align*}
We are going to prove that the above integrands satisfy the
hypotheses of Proposition \ref{p1.2}. Indeed, for any
$p\in[2,\infty]$, the rules of the Malliavin derivative, the
hypercontractivity inequality and $(\ref{2.2.4})$ yield the
following.

\begin{align}
&\left\|\MultipleIntegral^{(n'),H}_{t}
\left(h(\cdot,t)\right)-\MultipleIntegral^{(n'),H}_{s}
\left(h(\cdot,s)\right)\right\|_{\SobolevSpace^{1,p}}\notag\\
&\leq
\left\|\int_{\Intervalt^{n'}}\left(\khst^{\ast,(n')}h\left(\cdot,t\right)\right)
\left(\theta_1,\ldots,\theta_{n'}\right)dW_{\theta_1}\cdots
dW_{\theta_{n'}}\right\|_{\SobolevSpace^{1,p}}\notag\\
&\quad+\left\|\int_{\left[0,s\right]^{n'}}\left(K_{H,s}^{\ast,(n')}\left[\Delta^{n'+1}h\left(\cdot,s;t\right)\right]\right)
\left(\theta_1,\ldots,\theta_{n'}\right)dW_{\theta_1}\cdots
dW_{\theta_{n'}}\right\|_{\SobolevSpace^{1,p}}\notag\\
&\leq
C(n',p)\left[\left\|\khst^{\ast,(n')}h\left(\cdot,t\right)\right\|_{L^{2}\left(\Intervalt^{n'}\right)}+\left\|K_{H,s}^{\ast,(n')}
\left[\Delta^{n'+1}h\left(\cdot,s;t\right)\right]\right\|_{L^2\left(\left[0,s\right]^{n'}\right)}\right]\notag\\
&\leq
C(n',p,T,\lambda,H)\left[\left\|h\right\|_{\mathbb{H}^\lambda(\Lambda^{(n')}_t)}\left|t-s\right|^H
+\left\|h\right\|_{\mathbb{H}^\lambda(\Lambda^{(n')}_t)}\left|t-s\right|^\lambda\right] \notag\\
&\leq
C(n',p,T,\lambda,H)\left\|h\right\|_{\mathbb{H}^\lambda(\Lambda^{(n')}_t)}\left|t-s\right|^{H
\wedge \lambda}.\label{2.5.0}
\end{align}

Thus, the $\SobolevSpace^{1,p}$--valued stochastic process
$Y^{(n')}=\left\{\MultipleIntegral^{(n'),H}_{t}\left(h(\cdot,t)\right);\;t\in
\IntervalT\right\}$ has $(H \wedge \lambda)$--H\"older continuous
sample paths, a.s.  Notice that $(H \wedge \lambda)+H>\frac{1}{2}$. By applying Proposition \ref{p1.2}, we can
write
\begin{align*}
\MultipleIntegral^{(n'+1),H}_{t}\left(h\right)-\MultipleIntegral^{(n'+1),H}_{s}\left(h\right)
&=\int_0^t \left(K_{H,t}^{\ast,(1)}\MultipleIntegral^{(n'),H}_{\star}\left(h(\cdot,\star)\right)\right)(\theta_{n'+1})dW_{\theta_{n'+1}}\\
&-\int_0^s\left(K_{H,s}^{\ast,(1)}\MultipleIntegral^{(n'),H}_{\star}\left(h(\cdot,\star)\right)\right)(\theta_{n'+1})dW_{\theta_{n'+1}}.
\end{align*}
Owing to (\ref{u}) applied to $h(\theta):=
\MultipleIntegral^{(n'),H}_ {\theta}\left(h(\cdot,\theta)\right)$
and the identity (\ref{50}), the last term can be decomposed into
the sum of four terms denoted by $N_i$, $i=1,\dots,4$, and defined
as follows:
\begin{align*}
&N_1=\int_s^t \MultipleIntegral^{(n'),H}_{\theta_{n'+1}}
\left(h(\cdot,\theta_{n'+1})\right)\kh\left(t,\theta_{n'+1}\right)dW_{\theta_{n'+1}},\\
&N_2=\int_s^t\Big(\int_{\theta_{n'+1}}^t\left[\MultipleIntegral^{(n'),H}_{r_{n'+1}}
\left(h(\cdot,r_{n'+1})\right)-\MultipleIntegral^{(n'),H}_{\theta_{n'+1}}
\left(h(\cdot,r_{n'+1})\right)\right] K_H(dr_{n'+1},\theta_{n'+1})\Big)dW_{\theta_{n'+1}},\\
&N_3=\int_s^t\Big(\int_{\theta_{n'+1}}^t\MultipleIntegral^{(n'),H}_{\theta_{n'+1}}
\left(\Delta^{n'+1}
h(\cdot,\theta_{n'+1};r_{n'+1})\right) K_H(dr_{n'+1},\theta_{n'+1})\Big)dW_{\theta_{n'+1}},\\
&N_4=\int_0^s \Big(\int_s^t\MultipleIntegral^{(n'),H}_{r_{n'+1}}
\left(h(\cdot,r_{n'+1})\right)K_H(dr_{n'+1},\theta_{n'+1})\;\Ind_{\left\{s\neq0\right\}}\Big)dW_{\theta_{n'+1}}.
\end{align*}

We now apply the induction assumption to write all these terms by means of $(n'+1)$--multiple integrals,
obtaining
\begin{align}
&N_1=\int_s^t
\Big(\int_{\left[0,\theta_{n'+1}\right]^{n'}}\left(K_{H,\theta_{n'+1}}^{\ast,(n')}
h\left(\cdot,\theta_{n'+1}\right)\right)\left(\theta_1,\ldots,\theta_{n'}\right)dW_{\theta_1}\cdots
dW_{\theta_{n'}}\Big)\nonumber\\
&\quad\quad\times\kh\left(t,\theta_{n'+1}\right)dW_{\theta_{n'+1}},\nonumber\\
&N_2=\int_s^t\Big(\int_{\theta_{n'+1}}^t\Big(\int_{\left[0,r_{n'+1}\right]^{n'}}
\left(K_{H,\theta_{n'+1},r_{n'+1}}^{\ast,(n')}
h\left(\cdot,r_{n'+1}\right)\right)\left(\theta_1,\ldots,\theta_{n'}\right)dW_{\theta_1}\cdots
dW_{\theta_{n'}}\Big)\nonumber\\
&\quad \quad \times K_H(dr_{n'+1},\theta_{n'+1})\Big)dW_{\theta_{n'+1}},\label{2.5.1}\\
&N_3=\int_s^t\Big(\int_{\theta_{n'+1}}^t\Big(\int_{\left[0,\theta_{n'+1}\right]^{n'}}
\left(K_{H,\theta_{n'+1}}^{\ast,(n')} \left[\Delta^{n'+1}
h(\cdot,\theta_{n'+1};r_{n'+1})\right]\right)\left(\theta_1,\ldots,\theta_{n'}\right)\nonumber\\
&\quad \quad \times dW_{\theta_1}\cdots dW_{\theta_{n'}}
\Big)K_H(dr_{n'+1},\theta_{n'+1})\Big)dW_{\theta_{n'+1}},\label{2.5.3}\\
&N_4=\int_0^s
\Big(\int_s^t\Big(\int_{\left[0,r_{n'+1}\right]^{n'}}\left(K_{H,r_{n'+1}}^{\ast,(n')}
h\left(\cdot,r_{n'+1}\right)\right)\left(\theta_1,\ldots,\theta_{n'}\right)dW_{\theta_1}\cdots
dW_{\theta_{n'}} \Big)\nonumber\\
&\quad \quad \times
K_H(dr_{n'+1},\theta_{n'+1})\;\Ind_{\left\{s\neq0\right\}}\Big)dW_{\theta_{n'+1}}.\label{2.5.2}
\end{align}
Finally, by the stochastic Fubini theorem applied to $(\ref{2.5.1})$,
$(\ref{2.5.2})$ and $(\ref{2.5.3})$, we obtain
\begin{align*}
&\MultipleIntegral^{(n'+1),H}_{t}
\left(h\right)-\MultipleIntegral^{(n'+1),H}_{s}
\left(h\right)\\
&=\int_{\left[0,\theta_{n'+1}\right]^{n'}\times
\left]s,t\right[}\left(K_{H,\theta_{n'+1}}^{\ast,(n')}
h\left(\cdot,\theta_{n'+1}\right)\right)\left(\theta_1,\ldots,\theta_{n'}\right)\kh\left(t,\theta_{n'+1}\right)\\
&\quad \times dW_{\theta_1}\cdots
dW_{\theta_{n'}}dW_{\theta_{n'+1}}\\
&\quad+\int_{\left[0,t\right]^{n'}\times\left]s,t\right[}\Big(\int_{\vee_{i=1}^{n'+1}\theta_i}^t
\left(K_{H,\theta_{n'+1},r_{n'+1}}^{\ast,(n')}
h\left(\cdot,r_{n'+1}\right)\right)\left(\theta_1,\ldots,\theta_{n'}\right)\\
&\quad\quad\times K_H(dr_{n'+1},\theta_{n'+1})\Big)dW_{\theta_1}\cdots
dW_{\theta_{n'}}dW_{\theta_{n'+1}}\\
&\quad+\int_{\left[0,\theta_{n'+1}\right]^{n'}\times\left]s,t\right[}\Big(\int_{\theta_{n'+1}}^t
\left(K_{H,\theta_{n'+1}}^{\ast,(n')} \left[\Delta^{n'+1}
h(\cdot,\theta_{n'+1};r_{n'+1})\right]\right)\left(\theta_1,\ldots,\theta_{n'}\right)\\
&\quad \quad \times
K_H(dr_{n'+1},\theta_{n'+1})\Big)dW_{\theta_1}\cdots
dW_{\theta_{n'}}dW_{\theta_{n'+1}}\\
&\quad+\int_{\left[0,t\right]^{n'}\times\left]0,s\right[}
\Big(\int_{\vee_{i=1}^{n'}\theta_i\vee
s}^t\left(K_{H,r_{n'+1}}^{\ast,(n')}
h\left(\cdot,r_{n'+1}\right)\right)\left(\theta_1,\ldots,\theta_{n'}\right)\\
&\quad \quad \times
K_H(dr_{n'+1},\theta_{n'+1})\;\Ind_{\left\{s\neq0\right\}}\Big)dW_{\theta_1}\cdots
dW_{\theta_{n'}}dW_{\theta_{n'+1}}.
\end{align*}
That is,
\begin{equation*}
\MultipleIntegral^{(n'+1),H}_{t}
\left(h\right)-\MultipleIntegral^{(n'+1),H}_{s}
\left(h\right)
=\int_{\Intervalt^{n'+1}}\left(\khst^{\ast,(n'+1)}h\right)
\left(\theta_1,\ldots,\theta_{n'+1}\right)dW_{\theta_1}\cdots
dW_{\theta_{n'+1}}
\end{equation*}
(see (\ref{2.2.2})). This ends the proof of (\ref{2.5}).

The upper bound estimate (\ref{2.5.0.0}) follows from (\ref{2.5}),
the hypercontractivity inequality and (\ref{2.2.4}), while the
conclusion on sample path regularity is a consequence of the
Kolmogorov's criterion. \hfill\qed

\medskip

 This is the analogue of Corollary \ref{c2.1} for  $H\in]\frac{1}{4},\frac{1}{2}[$.
 \begin{cor}
 \label{c2.2} Let $H\in]\frac{1}{4},\frac{1}{2}[$, $\lambda\in]0,1[$ be such that $\lambda+H>\frac{1}{2}$ and $\beta\in ]0,1[$. Let $h$ be a measurable function defined on $[0,T]^{n+1}$ such that
 the mapping $t\longmapsto h(\cdot,t)$
belongs to $\mathcal{C}^\beta\left([0,T];\nSingularSpace\right)$.
Then, the integral process
$\left\{\fnint_t\left(h(\cdot,t)\right),\;t \in\IntervalT\right\}$
given in Theorem $\ref{t2.2}$ has  $\gamma$--H\"older continuous
sample paths, a.s., with $\gamma\in]0,\beta \wedge H[$.
\end{cor}
{\it Proof}: It follows from the estimate (\ref{2.5.0}), with
$\lambda$ replaced now by $\beta$,  and Kolmogorov's criterion.
\hfill\qed

\section {Large deviation principle for Banach valued multiple integrals}
\label{s3} In this section we fix an integer $n\ge 1$ and we
consider the integral processes
$\left\{\fnint_t\left(h(\cdot,t)\right), t\in [0,T]\right\}$ given
in the Corollaries \ref{c2.1} and \ref{c2.2}, respectively, for a
fixed deterministic function $h$. The fractional Brownian motion
$B^H$ is replaced  by  $\sep \bh$, $\ep\in]0,\infty[$. Since the
integrand $h$ is deterministic, we obtain a family
$\left\{\ep^{\frac{n}{2}}\fnint_t\left(h(\cdot,t)\right), t\in
[0,T]\right\}$ indexed by $\ep\in]0,\infty[$. It can also be
thought of as a family of random variables taking values on the
space $\mathcal{C}^\gamma([0,T])$, for some value of
$\gamma\in]0,1[$ made explicit in the above mentioned corollaries.
Throughout this section we shall denote by $\left\{
\mathbf{I}^{(n),H,\ep}_h, \ep>0\right\}$ this object. Our purpose
is to establish a large deviation principle for the corresponding
family of probability laws. For this, we shall apply the results
of \cite{Le90}.

The first step consists of  describing the abstract Wiener space
to be considered in our problem. First, we should replace the
space of H\"older continuous functions by a Polish space by a
classical procedure (see for instance \cite{Ci60}). More
precisely, fix $\gamma\in]0,1[$ and for every $0<\delta\leq T$ set
\begin{equation*}
\omega_{x}(\delta)=\underset{s\neq t}{\underset{0\leq s, t \leq
T}{\sup_{\left|t-s\right|\leq\delta}}}
\frac{\left|x(t)-x(s)\right|}{\left|t-s\right|^\gamma}.
\end{equation*}

Let $\mathcal{C}^{\gamma,0}\left(\IntervalT\right)$
be the subspace of $\mathcal{C}^{\gamma}\left(\IntervalT\right)$ consisting of
functions such that $\lim_{\delta\rightarrow0^+}\omega_{x}(\delta)=0$. The space
$\mathcal{C}^{\gamma,0}\left(\IntervalT\right)$ is a closed
subspace of $\mathcal{C}^{\gamma}\left(\IntervalT\right)$, so that
endowed with the norm $\left\|\cdot\right\|_{\gamma}$ it is a separable
Banach space, whereas $\mathcal{C}^{\gamma}\left(\IntervalT\right)$ is not.

It is easy to check that $\mathcal{C}^\gamma([0,T])\subset \mathcal{C}^{\gamma',0}([0,T])$, for any $\gamma'<\gamma$.

Consider the set
\begin{equation}
\label{3.1}
\mathcal{H}^H = \left\{\varphi:\IntervalT\rightarrow
\mathbb{R}: \;\varphi(t)=\int_0^t\kh(t,\theta)\dot{\varphi}(\theta)d\theta, \dot{\varphi}\in L^2(\IntervalT) \right\},
\end{equation}
with $K_H$ given in (\ref{1.2}), endowed with the inner product
$
\left \langle \varphi_1,\varphi_2\right \rangle_{\mathcal{H}^H}=\left \langle
\dot{\varphi}_1,\dot{\varphi}_2\right \rangle_{L^2(\IntervalT)}.
$
This is the reproducing kernel Hilbert space of the fractional Brownian motion,
as it is proved in \cite{DU198}.
\begin{rem}
\label{r3.1}
For any $H\in]\frac{1}{2},1[$, we have
\begin{equation}
L^2(\IntervalT)\subset L^{\frac{1}{H}}(\IntervalT)\subset \left|
\mathcal{H}^H\right| \subset \mathcal{H}^H.
\end{equation}
Indeed, the first inclusion is obvious, while the second one
follows from the estimate $(\ref{1.8})$, and the last one is
pointed out in \cite{PT01}. This implies that $\mathcal{H}^H$
contains the set of continuous functions defined on $[0,T]$.

For $H\in]0,\frac{1}{2}[$, $\lambda\in]0,1[$ such that
$\lambda+H>\frac{1}{2}$, as is pointed out in \cite{Nu03},
\begin{equation}
\mathcal{C}^\lambda(\IntervalT)\subset \mathcal{H}^H.
\end{equation}
As a consequence, $\mathcal{H}^H$ contains the set of Lipschitz continuous functions.
\end{rem}

Let $\mathcal{C}_0\left(\IntervalT\right)$ be the set of
continuous functions defined on $[0,T]$ vanishing at the origin,
$\bar
i:\mathcal{H}^H\hookrightarrow\mathcal{C}_0\left(\IntervalT\right)$
the canonical inclusion, and $\mathbf{P}^H$ be the law of the
fractional Brownian motion on
$\mathcal{C}_0\left(\IntervalT\right)$. The quadruple
$\left(\mathcal{C}_0\left(\IntervalT\right),\mathcal{H}^H,\bar
i,\mathbf{P}^H\right)$ is an abstract Wiener space (see
\cite{DU198} for the proof of these results). This result can be
strengthened, as is stated in the next Proposition.
\begin{prop}
\label{p3.1}For any $\zeta< H$, the set $\mathcal{H}^H$ is
included in $\mathcal{C}^{\zeta,0}([0,T])$. Moreover, denoting by
$i:\mathcal{H}^H\hookrightarrow \mathcal{C}^{\zeta,0}([0,T])$ the
canonical embedding, the quadruple
$$\left(\mathcal{C}^{\zeta,0}([0,T]),\mathcal{H}^H,i,\mathbf{P}^H\right),$$
 is an abstract Wiener space.
\end{prop}
{\it Proof}:  Indeed, let $\varphi\in \mathcal{H}^H$ and $0\le s<t\le T$. Schwarz's inequality yields
\begin{align*}
|\varphi(t)-\varphi(s)|&=\left\vert\int_0^t\left(K_H(t,r)-K_H(s,r)\right) \dot{\varphi}(r) dr\right\vert\\
&\le \Vert \varphi\Vert_{\mathcal{H}^H}\Vert
K_H(t,\cdot)-K_H(s,\cdot)\Vert_{L^2([0,T])}\le \Vert
\varphi\Vert_{\mathcal{H}^H}C|t-s|^H.
\end{align*}
Consequently,
$$\Vert \varphi\Vert_{\mathcal{C}^H([0,T])}\le C\Vert \varphi\Vert_{\mathcal{H}^H}.$$
This proves that the mapping $i$ is a continuous embedding.

Fix $x\in\mathcal{C}^{\zeta,0}([0,T])$ and consider the
approximating sequence, $\left\{x^{(n)}, n\in\mathbb{N}\right\}$,
consisting of linear interpolations on the dyadic numbers. That
is,
\begin{equation}
\label{3.2}
x^{(n)}(t)=\sum_{j=0}^{2^n-1}\left[x(t_j^n)+\frac{2^n}{T}\left(t-t_j^n\right)\left(x(t_{j+1}^n)-x(t_j^n)\right)\right]
\Ind\underset{\left[t_j^n,t_{j+1}^n\right]} {\left(t\right)},
\end{equation}
where $t_j^n=\frac{jT}{2^n}$, $j=0,\ldots,2^n$. Clearly, each $x^n$ is a Lipschitz function, consequently
$\left\{x^{(n)}, n\in\mathbb{N}\right\} \subset \mathcal{H}^H$ (see Remark \ref{r3.1}).

It is an easy exercise to check that this sequence  converges to
$x$ in the norm $\Vert\cdot\Vert_{\zeta}$. Therefore,
$i\left(\mathcal{H}^H\right)$ is dense in
$\mathcal{C}^{\zeta,0}([0,T])$.

Finally, since the trajectories of the fractional Brownian motion
are a.s. in $\mathcal{C}^{\zeta,0}([0,T])$, for any
$\zeta\in]0,H[$, we have
$(\mathbf{P}^H)^\ast\left(\mathcal{C}^{\zeta,0}([0,T])\right)=1$,
where  $(\mathbf{P}^H)^\ast$ denotes the exterior measure. This
leads to the conclusion (see Theorem 2.4 of \cite{BBK92}).
\hfill\qed

\smallskip

\begin{rem}
 \label{r3.2} An alternate approximation sequence by Lipschitz functions for $x\in\mathcal{C}^{\zeta,0}([0,T])$ in the norm $\Vert\cdot\Vert_{\zeta}$, is provided by
\begin{equation*}
x^{(n)}(t)=\frac{n}{T}\int_t^{t+\frac{T}{n}}x(\tau)d\tau-\frac{n}{T}\int_0^\frac{T}{n}x(\tau)d\tau,
\end{equation*}
with $x(t)=x(T)$ for $t\ge T$ (see \cite{MS61}, page 276).
\end{rem}
In the sequel, we consider as reference probability space the
triple $\left(\mathcal{C}^{\gamma,0}([0,T]), \mathcal {B},
\mathbf{P}^H\right)$, where $\mathcal {B}$ is the Borel
$\sigma$--field of $\mathcal{C}^{\gamma,0}([0,T])$. In particular,
the random variables $\left\{ \mathbf{I}^{(n),H,\ep}_h,
\ep>0\right\}$ and $\mathbf{I}^{(n),H}_h:= \mathbf{I}^{(n),H,1}_h$
are supposed to be defined in this probability space. Under the
hypotheses of Corollaries \ref{c2.1} and \ref{c2.2}, depending
whether $H\in]\frac{1}{2},1[$ or $H\in]\frac{1}{4},\frac{1}{2}[$,
they are $\mathcal{C}^{\gamma,0}([0,T])$--valued random variables
with $\gamma\in]0,\beta_0[$, where $\beta_0=\beta \wedge
\left(H-1/q\right)$ if $H\in ]\frac{1}{2},1[$, $qH>1$,  and
$\beta_0= \beta\wedge H$, if $H\in]\frac{1}{4},\frac{1}{2}[$,
respectively. We also denote by $\mathbf{I}^{(n),H}_h(B^H +
\varphi)$, for $\varphi \in \mathcal{H}^H$, the multiple
integral with respect to $\bh+\varphi$ instead of $\bh$.

\smallskip

The main theorem of \cite{Le90} (see page 4) applied to the
abstract Wiener space given in Proposition \ref{p3.1}, the
separable Banach space $B=\mathcal{C}^{\gamma,0}([0,T])$ and the
random variable $\mathbf{I}^{(n),H}_h$ in the
$\mathcal{C}^{\gamma,0}([0,T])$--valued Wiener chaos of degree
$n$, implies the large deviation principle for the family of laws
of $\left\{ \mathbf{I}^{(n),H,\ep}_h, \ep>0\right\}$ given here.
\begin{thm}
\label{t3.1} Each one of the two sets of assumptions
\begin{itemize}
\item[(i)] $H\in]\frac{1}{2},1[$, $qH>1$, $\beta\in ]0,1[$, $h$ is
a measurable function defined on $[0,T]^{n+1}$ such that the
mapping $t\longmapsto h(\cdot,t)$ belongs to
$\mathcal{C}^\beta\left([0,T];L^q(\Lambda_T^{(n)})\right)$,

\item [(ii)] $H\in]\frac{1}{4},\frac{1}{2}[$, $\lambda\in]0,1[$ is
such that $\lambda+H>\frac{1}{2}$,  $\beta\in ]0,1[$, $h$ is a
measurable function defined on $[0,T]^{n+1}$ such that the mapping
$t\longmapsto h(\cdot,t)$ belongs to
$\mathcal{C}^\beta\left([0,T];\nSingularSpace\right)$,
\end{itemize}
yields the following:

\medskip
For any closed set $F\subset \mathcal{C}^{\gamma,0}([0,T])$,
$$
\limsup_{\ep\to 0} \ep \log \mathbf{P}^H \left( \mathbf{I}^{(n),H,\ep}_h\in F\right) \le -\mathcal{I}(F).
$$
For any open set $G\subset \mathcal{C}^{\gamma,0}([0,T])$,
$$
\liminf_{\ep\to 0} \ep \log \mathbf{P}^H \left( \mathbf{I}^{(n),H,\ep}_h\in F\right) \ge -\mathcal{I}(G).
$$
Here $\gamma\in]0,\beta\wedge \left(H-1/q\right)[$ under
assumptions $(i)$, $\gamma\in]0,\beta\wedge H[$, under $(ii)$,
respectively, and the rate functional $\mathcal{I}$ is defined  by
$$
\mathcal{I}(\Phi)= \frac{\Vert \varphi\Vert^2_{\mathcal{H}^H}}{2}, \,  \Phi\in \mathcal{C}^{\gamma,0}([0,T]),
$$
if there exist $\varphi\in\mathcal{H}^H$ such that $\Phi =
\mathbb{E}\left(\mathbf{I}^{(n),H}_h(B^H + \varphi)\right)$,
and $\mathcal{I}(\Phi)= \infty$, otherwise.
\end{thm}

The next discussion completes the description of the rate functional of the large deviation principle.

The construction of multiple integrals with respect to the
fractional Brownian motion and their properties proved in the
preceding sections have an analogue when the fractional Brownian
motion is replaced by a deterministic function  $\varphi\in
\mathcal{H}^H$. More precisely, set
$$
\mathbf{J}_t^{(n),H}(h)= \int_{\Lambda_t^{(n)}} h(\theta_1,
\dots,\theta_n) \varphi(d\theta_1)\cdots\varphi(d\theta_n), \,
t\in [0,T],
$$
then we have the following.
\begin{enumerate}
\item Let $H\in]\frac{1}{2},1[$ and $h\in
L^{\frac{1}{H}}(\Lambda_T^{(n)})$. The function $t\rightarrow
\mathbf{J}_t^{(n),H}(h)$ is well defined and satisfies
$$
\sup_{t\in[0,T]} \left\vert\mathbf{J}_t^{(n),H}(h)\right\vert\le C
\Vert h\Vert_{L^{\frac{1}{H}}(\Lambda_T^{(n)})} \Vert
\varphi\Vert_{\mathcal{H}^H}^n.
$$
\item Let $H\in]\frac{1}{4},\frac{1}{2}[$, $\lambda\in ]0,1[$ such
that $\lambda+H>\frac{1}{2}$, and $h\in
\mathbb{H}^\lambda(\Lambda_T^{(n)})$. The function $t\rightarrow
\mathbf{J}_t^{(n),H}(h)$ is well defined and satisfies
$$
\sup_{t\in[0,T]} \left\vert\mathbf{J}_t^{(n),H}(h)\right\vert\leq
C \Vert h\Vert_{\mathbb{H}^\lambda(\Lambda_T^{(n)})} \Vert
\varphi\Vert_{\mathcal{H}^H}^n.
$$
\end{enumerate}
Similarly, one can state the analogues of Corollaries \ref{c2.1} and \ref{c2.2} and prove the existence of $\mathbf{J}_h^{(n),H}$, that is, the function
$t\rightarrow \mathbf{J}_t^{(n),H}(h(\cdot,t))$ of $\mathcal{C}^{\gamma,0}([0,T])$, for  specific values of $\gamma$.

For the proof of these facts, we first establish recursively -as for the stochastic integrals- the representation
$$
\mathbf{J}_t^{(n),H}(h)= \int_{[0,t]^n} \left(K^{\ast,(n)}_{H,t}
h\right) (\theta_1,\dots,\theta_n) \dot\varphi(\theta_1)
\cdots\dot\varphi(\theta_n) d\theta_1 \cdots d\theta_n.
$$
Then, we apply Cauchy-Schwarz inequality and finally,
(\ref{2.1.0}) and (\ref{2.2.4}) (with s:=0), respectively.

The stochastic integral $\mathbf{I}_h^{(1),H}(B^H+\varphi)$
satisfies
$\mathbb{E}\left(\mathbf{I}_h^{(1),H}(B^H+\varphi)\right)=
\mathbf{J}_h^{(1),H}$. Applying this fact recursively yields
$$
 \mathbb{E}\left(\mathbf{I}^{(n),H}_h(B^H + \varphi)\right) = \mathbf{J}_h^{(n),H}.
$$
Hence, we recover the usual description of the rate functional in terms of the {\it skeleton} of the Gaussian functional.




\begin{thebibliography}{Bi1-81}

\bibitem{AMN01} Al\`os, E., Mazet, O., Nualart, D.:
\emph{Stochastic calculus with respect to Gaussian processes.}
Ann. Probab. \textbf{29} (2) (2001), 766-801.

\bibitem{BBK92} Baldi, P., Ben Arous, G., Kerkyacharian, G.:
\emph{Large deviations and the {S}trassen theorem in {H}\"older
norm.} Stochastic Processes and their Applications. \textbf{42}
(1) (1992), 171-180.

\bibitem{Ci60} Ciesielski, Z.: \emph{On the Isomorphisms of the Spaces $H_\alpha$ and
$m$.} Bull. Acad. Pol. Sci. \text{7} (4), (1960), 217-222.

\bibitem{D1} Decreusefond, L.: \emph{Regularity Properties of Some Stochastic Volterra Integrals
with Singular Kernel.} Potential Anal. \textbf{16}, (2002),
139-149.

\bibitem{D2} Decreusefond, L.: \emph{Stochastic integration with respect to Volterra processes.}
Ann. I. H. Poincar\'e - PR \textbf{41} (2005), 123-149.

\bibitem{DU198} Decreusefond, L., \"Ust\"unel, A. S.:
\emph{Stochastic analysis of the fractional Brownian motion.}
Potential Anal. \textbf{10}, (1998), 177-214.


\bibitem{Le90} Ledoux, M.: \emph{A note on large deviations for Wiener
chaos.} In: S\'eminaire de probabilit\'es, XXIV, Lecture Notes in
Math., \textbf{1426}. Springer Berlin-Heidelber-New York, (1990),
1-14.

\bibitem{LQZ} Ledoux, M., Qian, Z., Zhang, T.: \emph{Large deviations and support theorem for diffusion
processes via rough paths. } Stochastic Processes and their
applications \textbf{102} (2002), 265-283.

\bibitem{LT91} Ledoux, M., Talagrand, M.: \emph{Probability in Banach Spaces.}
Isoperimetry and processes. Ergebnisse der Mathematik und ihrer
Grenzgebiete (3) [Results in Mathematics and Related Areas (3)],
23. Springer-Verlag, (1991).

\bibitem{L-Q} Lyons, T. J., Qian, Z.: \emph{System control and
rough paths.} Oxford Mathematical Monographs. Oxford Science
Publications. Oxford University Press, Oxford (2002).

\bibitem{MWNPA92} Mayer-Wolf, E., Nualart, D.,
P\'{e}rez-Abreu, V.: \emph{Large deviations for multiple
Wiener-It\^o integral processes.} In: S\'eminaire de
probabilit\'es, XXVI, Lecture Notes in Math., \textbf{1526}.
Springer, (1992), 11-31.

\bibitem{millet-ss} Millet, A., Sanz-Sol\'e, M.: \emph{Large deviations for rough paths
of the fractional Brownian motion.} Ann. I. H. Poincar\'e - PR
\textbf{42} (2006), 245-271.

\bibitem{MS61} Musielak, J., Semadeni, Z.: \emph{Some classes of Banach spaces depending
on a parameter.} Studia Mathematica, T. XX, (1981), 271-284.

\bibitem{Nu03} Nualart, D.: \emph{Stochastic integration with respect to fractional Brownian
motion and applications.} In: Stochastic models. Contemp. Math.
\textbf{336}. Amer. Math. Soc. (2003), 3-39.

\bibitem{Nu06} Nualart, D.: \emph{The Malliavin Calculus and Related
Topics.} Probability and its Applications. Springer-Verlag, 2nd
Edition, (2006).

\bibitem{PAT02} P\'{e}rez-Abreu, V., Tudor, C.: \emph{Multiple Stochastic Fractional Integrals:
A Transfer Principle for Multiple Stochastic Fractional
Integrals.} Bol. Soc. Mat. Mexicana \textbf{8} (3) (2002),
187-203.

\bibitem{PT01} Pipiras V., Taqqu, M. S.: \emph{Are classes of deterministic
integrands for fractional Brownian motion on a interval complete?}
Bernoulli \textbf{7} (2001), 873-897.

\bibitem{SKM93} Samko, S. G., Kilbas, A. A., Marichev, O.
I.: \emph{Fractional Integrals and Derivatives.} Gordon and
Breach, New York. (1993).

\bibitem{Stein} Stein, E. M.: \emph{Harmonic Analysis: real-variable methods, orthogonality, and
oscillatory integrals.} Princeton Mathematical Series,
\textbf{43}. Princeton University Press, 2nd printing, (1995).

\end{thebibliography}
\end{document}